\documentclass[11pt]{amsart}
\usepackage[english]{babel}
\usepackage{amsthm} 
\usepackage{amssymb}
\usepackage{amsmath}
\usepackage{hyperref}
\usepackage{amsfonts}
\usepackage{verbatim}
\usepackage{stmaryrd}
\usepackage{fancyhdr}
\usepackage{mathrsfs}
\usepackage{graphicx}
\usepackage{color}
\usepackage{pdfsync}
\newtheorem{theorem}{Theorem}[section]

\newtheorem{definition}[theorem]{Definition}
\newtheorem{lemma}[theorem]{Lemma}
\newtheorem{proposition}[theorem]{Proposition}
\newtheorem{remark}{Remark}  

\numberwithin{equation}{section}

\newcommand{\N}{\mathbb{N}}
\newcommand{\R}{\mathbb{R}}

\newcommand{\rr}{\mathbb{R}}

\newcommand{\CC}{\mathbb{C}}
\newcommand{\cc}{\mathbb{C}}
\newcommand{\eps}{\varepsilon}

\newcommand {\be}{\begin{equation}}
\newcommand {\ee}{\end{equation}}
\newcommand {\ba}{\begin{array}}
\newcommand {\ea}{\end{array}}

\usepackage{varioref}
 
\def\Rom#1{\uppercase\expandafter{\romannumeral #1}}

\def\polhk#1{\setbox0=\hbox{#1}{\ooalign{\hidewidth \lower1.5ex\hbox{`}\hidewidth\crcr\unhbox0}}}

\begin{document}

\title[Generalized Mehler formula for time-dependent  quadratic operators]{Generalized Mehler formula for time-dependent non-selfadjoint quadratic operators and propagation of singularities}

\author{Karel Pravda-Starov}

\address{\noindent \textsc{Karel Pravda-Starov, IRMAR, CNRS UMR 6625, Universit\'e de Rennes 1, Campus de Beaulieu, 263 avenue du G\'en\'eral Leclerc, CS 74205,
35042 Rennes cedex, France}}
\email{karel.pravda-starov@univ-rennes1.fr}


\begin{abstract}
We study evolution equations associated to time-dependent dissipative non-selfadjoint quadratic operators. We prove that the solution operators to these non-autonomous evolution equations are given by Fourier integral operators whose kernels are Gaussian tempered distributions associated to non-negative complex symplectic linear transformations, and we derive a generalized Mehler formula for their Weyl symbols. 
Some applications to the study of the propagation of Gabor singularities (characterizing the lack of Schwartz regularity) for the solutions to non-autonomous quadratic evolution equations are given. 
\end{abstract}

\maketitle

\section{Introduction}

\subsection{Mehler formula and quadratic Hamiltonians}

In his seminal work~\cite{historic}, Ferdinand Gustav Mehler established in 1866 the following celebrated formula, since then known as Mehler formula
$$\sum_{\alpha \in \mathbb{N}^n}\phi_{\alpha}(x)\phi_{\alpha}(y)\omega^{|\alpha|}=\frac{1}{\pi^{\frac{n}{2}}(1-\omega^2)^{\frac{n}{2}}}\exp\Big(-\frac{1+\omega^2}{2(1-\omega^2)}(x^2+y^2)+\frac{2\omega}{1-\omega^2}x \cdot y\Big),$$
holding for all $\omega \in \cc$, $|\omega|<1$ and $x,y \in \rr^n$, where $(\phi_{\alpha})_{\alpha \in \mathbb{N}^n}$ stands for the Hermite orthonormal basis, see also e.g.~\cite[p.~20]{combescure} (Theorem~1). This formula has played a major role in mathematical physics and more specifically in quantum mechanics for the study of Schr\"odinger equations associated to quadratic Hamiltonians. 
It allows in particular to derive explicit formulas for the kernel 
$$K_t(x,y)=\frac{1}{\big(2\pi \sinh(2t)\big)^{\frac{n}{2}}}\exp\Big(-\frac{1}{2\sinh(2t)}\big((x^2+y^2)\cosh(2t)-2x\cdot y\big)\Big),$$
with $(x,y) \in \rr^{2n}$, $t>0$, and the Weyl symbol
$$a_t(x,\xi)=\frac{1}{(\cosh(t))^n}\exp\big(-(\xi^2+x^2)\tanh(t)\big),$$
with $(x,\xi) \in \rr^{2n}$, $t>0$, of the contraction semigroup $(e^{-tH})_{t \geq 0}$ on $L^2(\rr^n)$ generated by the harmonic oscillator 
$$H=-\Delta_x+x^2, \quad x \in \rr^n.$$

There are many works concerning the quantum evolutions generated by quadratic Hamiltonians and exact formulas, see e.g.~\cite{unterberger1,unterberger2}. 
Quadratic Hamiltonians are actually very important in partial differential equations as they provide non trivial examples of wave propagation phenomena, and in quantum mechanics. They also play a major role when studying the propagation of coherent states for general classes of real-valued Hamiltonians including Schr\"odinger operators with general potentials
$$-h^2\Delta_x+V(x),$$ 
as this propagation of coherent states can be approximated in the semi-classical limit by the quantum evolutions generated by time-dependent real-valued quadratic Hamiltonians, see e.g. the works by Combescure, Robert, Laptev and Sigal~\cite{robert1,combescure,laptev,robert2}. Indeed, time-dependent real-valued quadratic Hamiltonians naturally appear in these latter works as the Taylor expansion up to order two of general Hamiltonians\footnote{Even in the case when Hamiltonians actually do not depend on time} $H$
\begin{multline*}
\widehat{H}_{2}(t)=H(X(t))+(x-x(t)) \cdot \frac{\partial H}{\partial x}(X(t))+(D_x-\xi(t)) \cdot \frac{\partial H}{\partial \xi}(X(t))\\
+\frac{1}{2}(x-x(t),D_x-\xi(t))\Big(\frac{\partial^2H}{\partial X^2}(X(t))\Big)(x-x(t),D_x-\xi(t))^T,
\end{multline*}
with $D_x=i^{-1}\partial_x$, around the classical flows $X(t)=(x(t),\xi(t))$ given by Hamilton's equations
$$\dot x(t)=\frac{\partial H}{\partial \xi}(x(t),\xi(t)), \quad \dot \xi(t)=-\frac{\partial H}{\partial x}(x(t),\xi(t)).$$ 
Among many others, as for instance the understanding of the smoothing properties of quadratic evolution equations developed as an application in the present work, the above consideration is one important motivation for studying Schr\"odinger evolutions associated to time-dependent quadratic Hamiltonians.

In the self-adjoint case, that is, for Schr\"odinger equations associated to real-valued time-dependent quadratic Hamiltonians, the propagation of coherent states is now fully understood thanks to the works of Combescure, Robert and Hagedorn~\cite{combescure2,robert1,hage}. We also refer the readers to the recent book by Combescure and Robert~\cite{combescure} for a comprehensive overview on this topic and others references herein. The properties and the structure of the Schr\"odinger solution operators generated by time-dependent real-valued quadratic Hamiltonians are also now fully understood thanks to the remarkable formula for their Weyl symbols derived by Mehlig and Wilkinson in~\cite{mehlig}, and proved independently by different approaches by Combescure and Robert~\cite{cubo}, and de Gosson~\cite{gosson}. The Mehlig-Wilkinson formula is recalled in the next section. 

On the other hand, H\"ormander studies in the work~\cite{mehler} the Schr\"odinger solution operators generated by complex-valued quadratic Hamiltonians giving rise to non-selfadjoint quadratic operators in the case when Hamiltonians do not depend on the time variable. In this beautiful work, H\"ormander establishes a very general Mehler formula for the Weyl symbols of these solution operators in the non-selfadjoint case that will be recalled below. This generalized Mehler formula derived by H\"ormander is now a keystone in numerous problems in mathematics and mathematical physics as it allows to perform exact computations for many problems.

In the present work, we bridge the gap between these two series of works by extending the general Mehler formula derived by H\"ormander for non-selfadjoint quadratic operators to the non-autonomous case, when complex-valued quadratic Hamiltonians are allowed to depend on the time variable. We believe that the generalized Mehler formula derived in this paper will also become a cornerstone in coming works on non-autonomous general non-selfadjoint evolution equations as it is already the case in particular for the study of propagation of coherent states in the selfadjoint case. Some applications to the study of the propagation of Gabor singularities (characterizing the lack of Schwartz regularity) for the solutions to non-autonomous quadratic evolution equations are given in the second part of the article.

 \subsection{Quadratic operators}

We consider quadratic operators. This class of operators stands for pseudodifferential operators
\begin{equation}\label{3.1}
q^w(x,D_x)u(x) =\frac{1}{(2\pi)^n}\int_{\R^{2n}}{e^{i(x-y) \cdot \xi}q\Big(\frac{x+y}{2},\xi\Big)u(y)dyd\xi}, \quad n \geq 1, 
\end{equation}
defined by the Weyl quantization of complex-valued quadratic symbols 
\begin{align*}
q: \rr^{2n} &\rightarrow \cc, \\
(x,\xi) &\mapsto q(x,\xi).
\end{align*}
These non-selfadjoint operators are only differential operators since the Weyl quantization of the quadratic symbol
$x^{\alpha} \xi^{\beta}$, with $(\alpha,\beta) \in \N^{2n}$, $|\alpha+\beta|=2$, is simply given by
\begin{equation}\label{forsym}
(x^{\alpha} \xi^{\beta})^w=\textrm{Op}^w(x^{\alpha} \xi^{\beta})=\frac{x^{\alpha}D_x^{\beta}+D_x^{\beta} x^{\alpha}}{2}, 
\end{equation}
with $D_x=i^{-1}\partial_x$. The maximal closed realization of a quadratic operator $q^w(x,D_x)$ on $L^2(\rr^n)$, that is, the operator equipped with the domain
\begin{equation}\label{dom1}
D(q^w)=\big\{u \in L^2(\rr^n) : \ q^w(x,D_x)u \in L^2(\rr^n)\big\},
\end{equation}
where $q^w(x,D_x)u$ is defined in the distribution sense, is known to coincide with the graph closure of its restriction to the Schwartz space~\cite{mehler} (pp.~425-426),
$$q^w(x,D_x) : \mathscr{S}(\rr^n) \rightarrow \mathscr{S}(\rr^n).$$

When the real part of the symbol is non-positive $\textrm{Re }q \leq 0$, the quadratic operator $q^w(x,D_x)$ equipped with the domain (\ref{dom1}) is maximal dissipative and generates a strongly continuous contraction semigroup $(e^{t q^w})_{t \geq 0}$ on $L^2(\rr^n)$~\cite{mehler} (pp. 425-426). The classical theory of strongly continuous semigroups~\cite[Chapter~4]{pazy} then shows that the function 
$$u \in C^0([0,+\infty[,L^2(\rr^n)) \cap C^1(]0,+\infty[,L^2(\rr^n)),$$ 
defined by $u(t)=e^{tq^w}u_0$ when $t \geq 0$, with $u_0 \in D(q^w)$, satisfies $u(t) \in D(q^w)$ for all $t \geq 0$, and is a classical solution to the autonomous Cauchy problem 
\begin{equation}\label{ty0}
\left\lbrace\begin{array}{l}
\frac{du(t)}{dt}=q^w(x,D_x)u(t), \qquad  t \geq 0,\\
u(0)=u_0.
\end{array}\right.
\end{equation}
Furthermore, the solution operator $e^{t q^w}$ for $t \geq 0$, is shown in~\cite{mehler} (Theorem~5.12) to be a Fourier integral operator $\mathscr{K}_{e^{2itF}}$, whose kernel is a Gaussian tempered distribution $K_{e^{2itF}} \in \mathscr{S}'(\rr^{2n})$ associated to the non-negative complex symplectic linear transformation
$$e^{2itF} : \cc^{2n} \rightarrow \cc^{2n},$$
where $F$ denotes the Hamilton map of the quadratic form $q$. This Hamilton map is the unique matrix $F \in \CC^{2n \times 2n}$ satisfying the identity
\begin{equation}\label{10}
\forall  (x,\xi) \in \R^{2n},\forall (y,\eta) \in \R^{2n}, \quad q((x,\xi),(y,\eta))=\sigma((x,\xi),F(y,\eta)), 
\end{equation}
with $q(\cdot,\cdot)$ the polarized form associated to $q$, where $\sigma$ stands for the standard symplectic form
\begin{equation}\label{11}
\sigma((x,\xi),(y,\eta))=\langle \xi, y \rangle -\langle x, \eta\rangle=\sum_{j=1}^n(\xi_j y_j-x_j \eta_j),
\end{equation}
with $x=(x_1,...,x_n)$, $y=(y_1,....,y_n)$, $\xi=(\xi_1,...,\xi_n)$, $\eta=(\eta_1,...,\eta_n) \in \cc^n$. 
In this work, the notation
$$\langle x,y \rangle=\sum_{j=1}^nx_jy_j, \quad x=(x_1,...,x_n) \in \cc^n, \ y=(y_1,...,y_n) \in \cc^n,$$
denotes the inner product on $\cc^n$, which is linear in both variables and not sesquilinear. 
We notice that a Hamilton map is skew-symmetric with respect to the symplectic form
\begin{multline}\label{a1}
\sigma((x,\xi),F(y,\eta))=q((x,\xi),(y,\eta))=q((y,\eta),(x,\xi))\\
=\sigma((y,\eta),F(x,\xi))=-\sigma(F(x,\xi),(y,\eta)), 
\end{multline}
by symmetry of the polarized form and skew-symmetry of the symplectic form. The Hamilton map $F$ is given by
\begin{equation}\label{jk1}
F=\sigma Q, 
\end{equation}
if $Q \in \cc^{2n\times 2n}$ denotes the symmetric matrix defining the quadratic form $q(X)=\langle Q X,X\rangle$, with $X=(x,\xi) \in \rr^{2n}$, and 
$$\sigma = \left(\begin{array}{cc}
0 & I_n \\
-I_n & 0
\end{array}
\right) \in \rr^{2n \times 2n},$$
with $I_n \in \rr^{n \times n}$ the identity matrix. The definition and the basic properties of the class of Fourier integral operators $\mathscr{K}_{\mathcal{T}}$, whose kernels $K_{\mathcal{T}} \in \mathscr{S}'(\rr^{2n})$ are Gaussian tempered distributions associated to non-negative complex symplectic linear transformations $\mathcal{T}$ are given in Section~\ref{FIO}.

On the other hand, H\"ormander shows in~\cite{mehler} (Theorem~4.2) that the solution operator $e^{t q^w}$ for $t \geq 0$, can also be considered as a pseudodifferential operator defined by the Weyl quantization of a tempered symbol $p_t \in \mathscr{S}'(\rr^{2n})$ explicitly given by the celebrated general Mehler formula
\begin{equation}\label{ty5}
p_t(X)=\frac{1}{\sqrt{\det(\cos tF)}}e^{\sigma(X,\tan(tF)X)} \in L^{\infty}(\rr^{2n}), \qquad X=(x,\xi) \in \rr^{2n},
\end{equation}
whenever the time variable $t \geq 0$ obeys the condition $\det(\cos tF) \neq 0$. 
Under the sole assumption that the real part of the symbol is non-positive $\textrm{Re }q \leq 0$, this condition $\det(\cos tF) \neq 0$ is not always satisfied. According to~\cite[p. 427]{mehler}, it is for instance the case of the solution operator associated to the harmonic Schr\"odinger operator $(e^{-it(D_x^2+x^2)})_{t \in \rr}$, whose Weyl symbol is given by
$$(x,\xi) \mapsto \frac{1}{\cos t}e^{-i(\xi^2+x^2)\tan t} \in L^{\infty}(\rr^{2n}),$$ 
when $\cos t \neq 0$, whereas when $t=\frac{\pi}{2}+k\pi$ with $k \in \mathbb{Z}$, it is given by the Dirac mass
\begin{equation}\label{sx1}
i(-1)^{k+1}\pi \delta_0(x,\xi) \in \mathscr{S}'(\rr^{2n}).
\end{equation}
The above formula accounts in particular for phenomena of reconstruction of singularities known for the Schr\"odinger equation~\cite{weinstein,zelditch1,zelditch2}.

In the present work, we unveil how the general Mehler formula (\ref{ty5}) extends to the non-autonomous case.

\subsection{Statements of the main results}

We consider time-dependent quadratic operators $q_t^w(x,D_x)$ whose symbols have coefficients 
$$q_t(x,\xi)=\sum_{\substack{\alpha, \beta \in \mathbb{N}^n \\ |\alpha+\beta| = 2}}(q_t)_{\alpha,\beta}x^{\alpha} \xi^{\beta},$$
depending continuously on the time variable $0 \leq t \leq T$, with $T>0$, and non-positive real parts 
\begin{equation}\label{ty3}
\textrm{Re }q_t \leq 0, \qquad 0 \leq t \leq T.
\end{equation}
We study the non-autonomous Cauchy problem 
\begin{equation}\label{ty1}
\left\lbrace\begin{array}{l}
\frac{du(t)}{dt}=q_t^w(x,D_x)u(t), \qquad 0 < t \leq T,\\
u(0)=u_0.
\end{array}\right.
\end{equation}
A continuous function $u \in C^0([0,T],L^2(\rr^n))$ is a classical solution of (\ref{ty1}) if $u$ is continuously differentiable in $L^2(\rr^n)$ on $]0,T]$, verifies $u(t) \in D(q_t^w)$ for all $0<t \leq T$, and satisfies the Cauchy problem (\ref{ty1}) in $L^2(\rr^n)$. As mentioned in~\cite[p.~139]{pazy}, there are no simple conditions that guarantee the existence of classical solutions for abstract non-autonomous Cauchy problems as (\ref{ty1}). Following~\cite[Definition~5.4.1]{pazy}, we therefore restrict ourselves to the study of a restricted notion of solutions. Setting
\begin{equation}\label{hilbert}
B=\{u \in  L^2(\rr^n) : x^{\alpha}D_x^{\beta}u \in L^2(\rr^n), \ \alpha, \beta \in \mathbb{N}^n, \ |\alpha+\beta| \leq 2\},
\end{equation}
the Hilbert space equipped with the norm 
$$\|u\|_{B}^2=\sum_{\substack{\alpha, \beta \in \mathbb{N}^n \\ |\alpha+\beta| \leq 2}}\|x^{\alpha}D_x^{\beta}u\|_{L^2(\rr^n)}^2,$$
we consider the following notion of $B$-valued solutions:

\bigskip

\begin{definition}\label{def1} \emph{($B$-valued solutions)}.
A continuous function $u \in C^0([0,T],B)$ is a $B$-valued solution of the non-autonomous Cauchy problem (\ref{ty1}) if $u \in C^1(]0,T],L^2(\rr^n))$ and (\ref{ty1}) is satisfied in $L^2(\rr^n)$.
\end{definition}

\bigskip

A $B$-valued solution differs from a classical solution by satisfying $u(t) \in B \subset D(q_t^w)$ for all $0 \leq t \leq T$, rather than only $u(t) \in D(q_t^w)$, and by being continuous in the stronger $B$-norm rather than merely in the $L^2(\rr^n)$-norm.

The first result contained in this paper establishes the existence and uniqueness of $B$-valued solutions to the non-autonomous Cauchy problem (\ref{ty1}):

\bigskip

\begin{theorem}\label{th1} \emph{(Existence and uniqueness of $B$-valued solutions).}
Let $T>0$ and $q_t : \rr^{2n} \rightarrow \cc$ be a time-dependent complex-valued quadratic form with a non-positive real part $\emph{\textrm{Re }}q_t \leq 0$ for all 
$0 \leq t \leq T$, and
whose coefficients depend continuously on the time variable $0 \leq t \leq T$, then for every $u_0 \in B$, the non-autonomous Cauchy problem
$$\left\lbrace\begin{array}{l}
\frac{du(t)}{dt}=q_t^w(x,D_x)u(t), \qquad 0 < t \leq T,\\
u(0)=u_0,
\end{array}\right.$$
has a unique $B$-valued solution. This solution is given by $u(t)=\mathscr{U}(t,0)u_0$ for all $0 \leq t \leq T$,
where $(\mathscr{U}(t,\tau))_{0 \leq \tau \leq t \leq T}$ is a contraction evolution system on $L^2(\rr^n)$, that is, a two parameters family of bounded linear operators on $L^2(\rr^n)$ satisfying
\begin{itemize}
\item[$(i)$] $\mathscr{U}( \tau, \tau)=I_{L^2(\rr^n)}$,\ \ $\mathscr{U}(t,r)\mathscr{U}(r, \tau)=\mathscr{U}(t, \tau)$ for all $0 \leq  \tau \leq r \leq t \leq T$
\item[$(ii)$] $(t, \tau) \mapsto \mathscr{U}(t, \tau)$ is strongly continuous on $L^2(\rr^n)$ for all $0 \leq  \tau \leq t \leq T$
\item[$(iii)$] $\forall 0 \leq \tau \leq t \leq T, \ \|\mathscr{U}(t, \tau)\|_{\mathcal{L}(L^2)} \leq 1$, with $\|\cdot\|_{\mathcal{L}(L^2)}$ standing for the operator norm on $L^2(\rr^n)$
\end{itemize}
\end{theorem}

\bigskip

In the autonomous case, we recall from~\cite{mehler} (Theorem~5.12) that the solution operator $e^{tq^w}$ for $t \geq 0$, is a Fourier integral operator whose kernel is a Gaussian tempered distribution associated to the non-negative complex symplectic linear transformation 
$e^{2itF} : \cc^{2n} \rightarrow \cc^{2n},$
where $F$ denotes the Hamilton map of the quadratic symbol $q$. The following result extends this description to the non-autonomous case, and shows that the evolution operators $\mathscr{U}(t, \tau)$, with $0 \leq \tau \leq t \leq T$, given by Theorem~\ref{th1} are also Fourier integral operators whose kernels are anew Gaussian tempered distributions associated to non-negative complex symplectic linear transformations:

\bigskip 

\begin{theorem}\label{th2} \emph{(Evolution operators as Fourier integral operators).}
Under the assumptions of Theorem~\ref{th1}, the evolution operator 
$$\mathscr{U}(t, \tau)=\mathscr{K}_{R(t, \tau)} : L^2(\rr^n) \rightarrow L^2(\rr^n), \qquad 0 \leq \tau \leq t \leq T,$$ 
is a Fourier integral operator whose kernel $K_{R(t, \tau)} \in \mathscr{S}'(\rr^{2n})$ is the Gaussian tempered distribution defined in the sense of Proposition~\ref{FIOdef} (Section~\ref{FIO}) associated to the non-negative complex symplectic linear transformation $R(t, \tau)$ given by the resolvent
\begin{equation}\label{resolvent}
\left\lbrace 
\begin{array}{ll}
\frac{d}{dt}R(t, \tau)=2iF_tR(t, \tau), \qquad 0 \leq t \leq T,\\
R(\tau,\tau)=I_{2n},
\end{array}
\right.
\end{equation}
with $0 \leq \tau \leq T$, where $F_t$ denotes the Hamilton map of $q_t$ and $I_{2n}$ stands for the $2n \times 2n$ identity matrix. 
On the other hand, the adjoint of the evolution operator 
$$\mathscr{U}(t, \tau)^*=\mathscr{K}_{\overline{R(t, \tau)}^{-1}} : L^2(\rr^n) \rightarrow L^2(\rr^n), \qquad 0 \leq \tau \leq t \leq T,$$ 
is the Fourier integral operator whose kernel $K_{\overline{R(t, \tau)}^{-1}} \in \mathscr{S}'(\rr^{2n})$ is the Gaussian tempered distribution  associated to the non-negative complex symplectic linear transformation $\overline{R(t, \tau)}^{-1}$.
Furthermore, the evolution operator 
$$\mathscr{U}(t, \tau)=\mathscr{K}_{R(t, \tau)} : \mathscr{S}(\rr^n) \rightarrow \mathscr{S}(\rr^n), \qquad 0 \leq \tau \leq t \leq T,$$
defines a continuous mapping on the Schwartz space which can be extended by duality as a continuous mapping on the space of tempered distributions
$$\mathscr{U}(t, \tau) : \mathscr{S}'(\rr^n) \rightarrow \mathscr{S}'(\rr^n), \qquad 0 \leq \tau \leq t \leq T,$$
defined as
$$\forall u \in \mathscr{S}'(\rr^n), \forall v \in \mathscr{S}(\rr^n), \quad \langle \mathscr{U}(t, \tau)u,\overline{v} \rangle_{\mathscr{S}'(\rr^n),\mathscr{S}(\rr^n)}=\langle u,\overline{\mathscr{U}(t, \tau)^*v} \rangle_{\mathscr{S}'(\rr^n),\mathscr{S}(\rr^n)}.$$
\end{theorem}

\bigskip

This description of the evolution operators as Fourier integral operators plays a major role below for studying the propagation of Gabor singularities for $B$-valued solutions to non-autonomous Cauchy problems (\ref{ty1}). Before studying this problem of propagation of singularities, we establish that the celebrated Mehler formula (\ref{ty5}) can also be extended to the non-autonomous case:

\bigskip

\begin{theorem}\label{th3} \emph{(Generalized Mehler formula for time-dependent quadratic Hamiltonians).}
Under the assumptions of Theorem~\ref{th1}, there exists a positive constant $\delta>0$ such that for all $0 \leq \tau \leq t \leq T$ and $0 \leq t-\tau < \delta$, the evolution operator 
$$\mathscr{U}(t, \tau)=p_{t,\tau}^w(x,D_x) : L^2(\rr^n) \rightarrow L^2(\rr^n),$$ 
is a pseudodifferential operator whose Weyl symbol $p_{t,\tau}$ is a  $L^{\infty}(\rr^{2n})$-function given by
$$p_{t,\tau}(X)=\frac{2^n}{\sqrt{\emph{\textrm{det}}\big(R(t,\tau)+I_{2n}\big)}}\exp\big(-i\sigma(X,\big(R(t,\tau)-I_{2n}\big)\big(R(t,\tau)+I_{2n}\big)^{-1}X\big)\big),$$
with $X=(x,\xi) \in \rr^{2n}$, where $R(t,\tau)$ denotes the resolvent defined in (\ref{resolvent}), $\sqrt{z}=e^{\frac{1}{2}\log z}$ with $\log$ the principal determination of the complex logarithm on $\cc \setminus \rr_-$, and where the quadratic form
$$X=(x,\xi) \in \rr^{2n} \mapsto -i\sigma(X,\big(R(t,\tau)-I_{2n}\big)\big(R(t,\tau)+I_{2n}\big)^{-1}X\big) \in \cc,$$
has a non-positive real part for all $0 \leq \tau \leq t \leq T$, $0 \leq t-\tau < \delta$.
\end{theorem}

\bigskip

In the autonomous case, that is, when $F_t=F$ for all $0 \leq t \leq T$, Theorem~\ref{th3} allows to recover the classical Mehler formula (\ref{ty5}). In this case, the resolvent $R(t,0)$ is indeed equal to $e^{2itF}$, and we observe that 
\begin{multline*}
-i\big(R(t,0)-I_{2n}\big)\big(R(t,0)+I_{2n}\big)^{-1}=-i(e^{2itF}-I_{2n})(e^{2itF}+I_{2n})^{-1}\\
=\sin(tF)\cos(tF)^{-1}=\tan(tF)
\end{multline*}
and
\begin{multline*}
2^{-2n}\textrm{det}\big(R(t,0)+I_{2n}\big)=2^{-2n}\textrm{det}(e^{2itF}+I_{2n})=2^{-2n}\textrm{det}(2\cos(tF)e^{itF})\\
=\textrm{det}(\cos(tF))e^{it\textrm{Tr}F}=\textrm{det}(\cos(tF)),
\end{multline*}
since by (\ref{jk1}), the trace of a Hamilton map $F=\sigma Q$ is zero $\textrm{Tr}(F)=0$, because 
\begin{equation}\label{vb1}
\textrm{Tr}(F)=\textrm{Tr}(F^T)=\textrm{Tr}(\sigma Q)=\textrm{Tr}(Q^T\sigma^T)=-\textrm{Tr}(Q\sigma)=-\textrm{Tr}(\sigma Q),
\end{equation}
by symmetry and skew-symmetry of the matrices $Q=Q^T$ and $\sigma^T=-\sigma$. As in the autonomous case (\ref{ty5}), notice that the Weyl symbol of
the evolution operator $\mathscr{U}(t, \tau)$ is not necessarily a $L^{\infty}(\rr^{2n})$-function for all $0 \leq \tau \leq t \leq T$. It accounts for the condition $0 \leq t-\tau < \delta$ appearing in the statement of Theorem~\ref{th3} to ensure that the determinant $\textrm{det}(R(t,\tau)+I_{2n}) \neq 0$ is non-zero and its square root well-defined.     

Let us now explain how the result of Theorem~\ref{th3} relates to the remarkable formula derived by Mehlig and Wilkinson in~\cite{mehlig}, and proved independently by different approaches by Combescure and Robert~\cite{cubo}, and de Gosson~\cite{gosson}. The Mehlig-Wilkinson formula provides the following explicit formula for the Weyl symbol 
$$R_G(X)=\frac{2^ne^{i\pi \nu}}{\sqrt{|\det(G+I_{2n})|}}\exp\big(-i\sigma(X,(G-I_{2n})(G+I_{2n})^{-1}X\big)\big),$$
with $X=(x,\xi) \in \rr^{2n}$,
of a metaplectic operator $\widehat{R}(G)$ associated to a real symplectic linear transformation $G : \rr^{2n} \rightarrow \rr^{2n}$ satisfying $\det(G+I_{2n}) \neq 0$, where the parameter $\nu \in \mathbb{Z}$ is an integer if $\det(G+I_{2n})>0$, or an half-integer $\nu \in \mathbb{Z}+\frac{1}{2}$ if $\det(G+I_{2n})<0$. The integer or half-integer $\nu$ is explicitly computed by de Gosson in~\cite{gosson}, and depends in particular in a non-trivial manner on the Maslov index of the metaplectic operator $\widehat{R}(G)$.

Under the assumptions of Theorem~\ref{th1}, we consider the case when the quadratic symbol~$q_t$ has a zero real part
$$\forall 0 \leq t \leq T, \quad \textrm{Re }q_t=0,$$
that is, when it writes as $q_t=i\tilde{q}_t$, with $\tilde{q}_t$ a real-valued quadratic form whose coefficients depend continuously on the time variable $0 \leq t \leq T$. The resolvent defined in (\ref{resolvent}) is in this case a real symplectic linear transformation $R(t,\tau) : \rr^{2n} \rightarrow \rr^{2n}$, and the evolution operator $\mathscr{U}(t, \tau)=\mathscr{K}_{R(t, \tau)}$ given by the associated Fourier integral operator is then known to be~\cite[p.~447-448]{mehler} a metaplectic operator associated to the real symplectic linear transformation $R(t,\tau)$. This accounts for the fact that in this specific case, the generalized Mehler formula derived in Theorem~\ref{th3} reduces to the Mehlig-Wilkinson formula for $G=R(t,\tau)$, where the parameter $\nu$ is here equal to zero due to continuity properties of the symbol and the smallness condition imposed on the parameter $0 \leq t-\tau < \delta$ in the statement of Theorem~\ref{th3}.

\subsection{Propagation of Gabor singularities}
By using the above description of the evolution operators as Fourier integral operators, we aim next at studying the possible (or lack of) Schwartz regularity for the $B$-valued solutions to non-autonomous Cauchy problems (\ref{ty1}).

The lack of Schwartz regularity of a tempered distribution is characterized by its Gabor wave front set whose definition and basic properties are recalled in appendix (Section~\ref{appendix}). The Gabor wave front set (or Gabor singularities) was introduced by H\"ormander~\cite{Hormander1} and measures the directions in the phase space in which a tempered distribution does not behave like a Schwartz function. It is hence empty if and only if a distribution that is a priori tempered is in fact a Schwartz function. The Gabor wave front set thus measures global regularity in the sense of both smoothness and decay at infinity.

\subsubsection{General case}

In the autonomous case, this question of propagation of Gabor singularities for the solutions to evolution equations
\begin{equation}\label{cp}
\left\{
\begin{array}{l}
\frac{du(t)}{dt}= q^w(x,D_x) u (t) , \qquad t \geq 0, \\
u(0)  = u_0 \in L^2(\rr^n),  
\end{array}
\right.
\end{equation}
associated to any dissipative quadratic operator was adressed by Rodino, Wahlberg and the author in the recent work~\cite{wahlberg}. In this work, it is pointed out that only Gabor singularities of the initial datum $u_0 \in L^2(\rr^n)$ contained in the singular space $S$ of the quadratic symbol $q$, can propagate for positive times along the curves given by the flow $(e^{-tH_{\textrm{Im}q}})_{t \in \rr}$ of the Hamilton vector field 
$$H_{\textrm{Im}q}=\frac{\partial \textrm{Im } q}{\partial \xi} \cdot \frac{\partial}{\partial_x}-\frac{\partial \textrm{Im }q}{\partial x} \cdot \frac{\partial}{\partial_{\xi}},$$
associated to the opposite of the imaginary part of the symbol. On the other hand, the Gabor singularities of the initial datum outside the singular space are all smoothed out for any positive time. More specifically, the following microlocal inclusion of Gabor wave front sets is established in~\cite{wahlberg} (Theorem~6.2),
\begin{equation}\label{etoile}
\forall u_0 \in L^2(\rr^n), \forall  t>0, \quad WF(e^{tq^w}u_0) \subset e^{-tH_{\textrm{Im}q}}\big(WF(u_0) \cap S\big) \subset S.
\end{equation}
The notion of singular space was introduced by Hitrik and the author in~\cite{Hitrik1} by pointing out the existence of a particular vector subspace in the phase space $\rr^{2n}$, which is intrinsically associated to a quadratic symbol $q$, and defined as the following finite intersection of kernels
\begin{equation}\label{h1}
S=\Big(\bigcap_{j=0}^{2n-1} \textrm{Ker}\big[\textrm{Re } F(\textrm{Im }F)^j \big]\Big) \cap \rr^{2n} \subset \rr^{2n}, 
\end{equation}
where $\textrm{Re }F$ and $\textrm{Im }F$ stand for the real and imaginary parts of the Hamilton map $F$ associated to $q$,
$$\textrm{Re }F=\frac{1}{2}(F+\overline{F}), \quad \textrm{Im }F=\frac{1}{2i}(F-\overline{F}),$$
which are respectively the Hamilton maps of the quadratic forms $\textrm{Re }q$ and $\textrm{Im }q$.
As pointed out in \cite{Hitrik1,kps11,kps21,viola1}, the singular space is playing a basic role in understanding the spectral and hypoelliptic properties of non-elliptic quadratic operators, as well as the spectral and pseudospectral properties of certain classes of degenerate doubly characteristic pseudodifferential operators~\cite{kps3,kps4,viola}. In the case when the singular space is zero $S=\{0\}$, the microlocal inclusion (\ref{etoile}) implies that the semigroup $(e^{tq^w})_{t \geq 0}$ enjoys regularizing properties of Schwartz type 
$$\forall u_0 \in L^2(\rr^{n}), \forall t>0, \quad e^{tq^w}u_0 \in \mathscr{S}(\rr^{n}),$$
for any positive time. It holds for instance for some non-selfadjoint non-elliptic kinetic operators as the Kramers-Fokker-Planck operator 
$$K=-\Delta_v+\frac{v^2}{4}+v\cdot\partial_x-\nabla V(x)\cdot\partial_v, \quad (x,v) \in \rr^{2},$$
with a quadratic potential $V(x)=ax^2$, $a \in \mathbb{R} \setminus \{0\}$, some operators appearing in models of finite-dimensional Markovian approximation of the general Langevin equation, or in chains of oscillators coupled to heat baths~\cite[Section~4]{kps11}.

In order to derive a microlocal inclusion for the propagation of Gabor singularities in the non-autonomous case, we need to generalize this notion of singular space to the time-dependent case. We consider the following definition:

\bigskip

\begin{definition}\label{def2}
Let $t_1 \leq t_2$ and $q_t : \rr^{2n} \rightarrow \cc$ be a time-dependent complex-valued quadratic form 
whose coefficients depend continuously on the time variable $t_1 \leq t \leq t_2$. The time-dependent singular space associated to the family of quadratic forms $(q_t)_{t_1 \leq t \leq t_2}$ is defined as 
\begin{equation}\label{ku10}
S_{t_1,t_2}=\Big( \bigcap_{t_1 \leq \tau \leq t_2}\emph{\textrm{Ker}}\big(\emph{\textrm{Im }}R(\tau,t_2)\big)\Big) \cap \rr^{2n},
\end{equation}
where $\emph{\textrm{Im }}R(t,\tau)=\frac{1}{2i}(R(t,\tau)-\overline{R(t,\tau)})$ denotes the imaginary part of the resolvent $R(t,\tau)$ defined in (\ref{resolvent}) and associated to the Hamilton map $F_t$ of $q_t$.
\end{definition}

\bigskip

When $q_t : \rr^{2n} \rightarrow \cc$ is a time-dependent complex-valued quadratic form with a non-positive real part $\textrm{Re }q_t \leq 0$ for all $t_1 \leq t \leq t_2$, with $t_1<t_2$, this definition truly extends the one given in the autonomous case. Indeed, when the quadratic form does not depend on time, that is, when $q_t=q$ for all $t_1 \leq t \leq t_2$, with $t_1<t_2$, we first observe from (\ref{ku10}) that the time-dependent singular space reduces to 
$$S_{t_1,t_2}=\Big( \bigcap_{t_1 \leq \tau \leq t_2}\textrm{Ker}\big(\textrm{Im }e^{-2i(t_2-\tau)F}\big)\Big) \cap \rr^{2n},$$
if $F$ denotes the Hamilton map of $q$, and recall from the proof of Theorem~6.2 in~\cite[formula (6.11)]{wahlberg} that we have
$$S=\Big(\bigcap_{j=0}^{2n-1} \textrm{Ker}\big[\textrm{Re } F(\textrm{Im }F)^j \big]\Big) \cap \rr^{2n}=\Big( \bigcap_{t_1 \leq \tau \leq t_2}\textrm{Ker}\big(\textrm{Im }e^{-2i(t_2-\tau)F}\big)\Big) \cap \rr^{2n}.$$
On the other hand, we also recall from the proof of Theorem~6.2 in~\cite[formula (6.18)]{wahlberg} that 
$$\forall t \in \rr, \quad e^{-tH_{\textrm{Im}q}}S=e^{-2t \textrm{Im }F}S=S.$$
The microlocal inclusion (\ref{etoile}) can therefore be rephrased as
\begin{equation}\label{etoile2}
\forall u_0 \in L^2(\rr^n), \forall  t>0, \quad WF(e^{tq^w}u_0) \subset e^{-tH_{\textrm{Im}q}}\big(WF(u_0)\big) \cap S.
\end{equation}
This microlocal inclusion of Gabor wave front sets can be extended to the non-autonomous case as follows:

\bigskip

\begin{theorem}\label{th4}
Under the assumptions of Theorem~\ref{th1}, the Gabor wave front set of the unique $B$-valued solution $u(t)=\mathscr{U}(t,0)u_0$ to the non-autonomous Cauchy problem
$$\left\lbrace\begin{array}{l}
\frac{du(t)}{dt}=q_t^w(x,D_x)u(t), \qquad 0 < t \leq T,\\
u(0)=u_0,
\end{array}\right.$$
with $u_0 \in B$, satisfies the microlocal inclusion
\begin{equation}\label{ku30}
\forall 0 \leq t \leq T, \quad WF(u(t)) \subset  \big(\emph{\textrm{Re }}R(t,0)\big)\big(WF(u_0)\big)  \cap S_{0,t},
\end{equation}
where $S_{0,t}$ is the time-dependent singular space associated to the family of quadratic forms $(q_{\tau})_{0 \leq \tau \leq t}$ and where $\emph{\textrm{Re }}R(t,0)=\frac{1}{2}(R(t,0)+\overline{R(t,0)})$ is the real part of the resolvent defined in (\ref{resolvent}).
\end{theorem}

\bigskip

As a direct consequence of Theorem~\ref{th4}, we observe that if there exists a positive time $0<t_0 \leq T$ such that the time-dependent singular space 
 is zero 
$$S_{0,t_0}=\Big( \bigcap_{0 \leq \tau \leq t_0}\textrm{Ker}\big(\textrm{Im }R(\tau,t_0)\big)\Big) \cap \rr^{2n}=\{0\},$$
then the non-autonomous Cauchy problem 
$$\left\lbrace\begin{array}{l}
\frac{du(t)}{dt}=q_t^w(x,D_x)u(t), \qquad 0 < t \leq T,\\
u(0)=u_0,
\end{array}\right.$$
enjoys regularizing properties of Schwartz type for all time $t_0 \leq t \leq T$, 
$$\forall u_0 \in B, \forall t_0 \leq t \leq T, \quad u(t)=\mathscr{U}(t,0)u_0 \in \mathscr{S}(\rr^{n}).$$
Indeed, we first deduce from (\ref{ku30}) and (\ref{nm1}) that $u(t_0)=\mathscr{U}(t_0,0)u_0 \in \mathscr{S}(\rr^n)$, since $WF(u(t_0)) \subset \rr^{2n} \setminus \{0\}$. By noticing from Theorem~\ref{th2} that the operator 
$$\mathscr{U}(t,t_0)=\mathscr{K}_{R(t,t_0)} : \mathscr{S}(\rr^n) \rightarrow \mathscr{S}(\rr^n),$$
is continuous, we finally obtain from Theorem~\ref{th1} that 
$$\forall u_0 \in B, \forall t_0 \leq t \leq T, \quad u(t)=\mathscr{U}(t,t_0)\underbrace{\mathscr{U}(t_0,0)u_0}_{u(t_0)\in  \mathscr{S}(\rr^{n})} \in \mathscr{S}(\rr^{n}).$$
The result of Theorem~\ref{th4} points out that no matter is the initial datum $u_0 \in B$, the possible Gabor singularities of $u(t)$ the solution at time $0 \leq t \leq T$ are all localized in the time-dependent singular space $S_{0,t}$. Furthermore, the possible Gabor singularities of the solution at time $t$ can only come from Gabor singularities of the initial datum which have propagated by the mapping given by the real part of the resolvent $\textrm{Re }R(t,0)$.

\subsubsection{Metaplectic case}
The general result of Theorem~\ref{th4} can be readily sharpened in the case when the quadratic symbol $q_t$ has a zero real part
$$\forall 0 \leq t \leq T, \quad \textrm{Re }q_t=0.$$
As mentioned above, the evolution operator $\mathscr{U}(t,0)=\mathscr{K}_{R(t,0)}$ is then a metaplectic operator associated to the real symplectic linear transformation 
$$R(t,0)=\textrm{Re }R(t,0) : \rr^{2n} \rightarrow \rr^{2n}.$$  
According to Definition~\ref{def2}, the time-dependent singular space $S_{0,t}=\rr^{2n}$ is then equal to the whole phase space since
$$\forall 0 \leq \tau \leq t,\quad \textrm{Im }R(\tau,t)=0,$$
and the symplectic invariance of the Gabor wave front set (\ref{sympinv}) directly implies that the solution satisfies 
$$\forall 0 \leq t \leq T, \quad WF(u(t))=  \big(\textrm{Re }R(t,0)\big)\big(WF(u_0)\big).$$
This sharpens the result of Theorem~\ref{th4} and extends the one obtained in~\cite{wahlberg} in the autonomous case.

\subsubsection{Outline of the article}

The article is organized in the following manner. Section~\ref{FIO} is devoted to recall the definition and the basic properties of Fourier integral operators associated to non-negative complex symplectic linear transformations. Section~\ref{proofth} provides the proof of the main results contained in this work (Theorems~\ref{th1},~\ref{th2} and~\ref{th3}), whereas Section~\ref{gabor} is devoted to the proof of the result of propagation of Gabor singularities (Theorem~\ref{th4}). Section~\ref{appendix} is an appendix recalling the definition and basic properties of the Gabor wave front set of a tempered distribution.

\section{Fourier integral operators associated to non-negative complex symplectic linear transformations}\label{FIO}

This section is devoted to recall the definition and the basic properties of Fourier integral operators 
associated to non-negative complex symplectic linear transformations.

This class of operators is used in~\cite{mehler} by H\"ormander to describe the properties of strongly continuous contraction semigroups $(e^{tq^w})_{t \geq 0}$ generated by maximal dissipative quadratic operators $q^w(x,D_x)$. Theorem~\ref{th2} points out that it also allows to describe the properties of evolution operators $\mathscr{U}(t, \tau)$ solving the non-autonomous Cauchy problems (\ref{ty1}).

In order to recall the definition of these operators, we closely follow the introduction to Gaussian calculus given in~\cite{mehler} (Section~5). Let $0 \neq u \in \mathscr{D}'(\rr^n)$ and set
$$\mathscr{L}_u=\Big\{L(x,\xi)=\sum_{j=1}^na_j\xi_j+\sum_{j=1}^nb_jx_j : L^w(x,D_x)u=0 \Big\}.$$
We recall that a distribution $u$ is said to be Gaussian if every distribution $v \in \mathscr{D}'(\rr^n)$ satisfying $L^w(x,D_x)v=0$ for all $L \in \mathscr{L}_u$, is necessarily a multiple of $u$.

Let $\mathcal{T} : \cc^{2n} \rightarrow \cc^{2n}$ be a non-negative complex symplectic linear transformation, that is, an isomorphism of $\cc^{2n}$ satisfying
$$\forall X,Y \in \cc^{2n}, \ \sigma(\mathcal{T}X,\mathcal{T}Y)=\sigma(X,Y); \quad \forall X \in \cc^{2n}, \ i\big(\sigma(\overline{\mathcal{T}X},\mathcal{T}X)-\sigma(\overline{X},X)\big) \geq 0.$$
Associated to this non-negative symplectic linear transformation is its twisted graph
\begin{equation}\label{graph}
\lambda_{\mathcal{T}}=\{(\mathcal{T}X,X') : X \in \cc^{2n} \} \subset \cc^{2n} \times \cc^{2n},
\end{equation}
where $X'=(x,-\xi) \in \cc^{2n}$, if $X=(x,\xi) \in \cc^{2n}$, which defines a non-negative Lagrangian plane of $\cc^{2n} \times \cc^{2n}$ equipped with the symplectic form 
$$\sigma_1((z_1,z_2),(\zeta_1,\zeta_2))=\sigma(z_1,\zeta_1)+\sigma(z_2,\zeta_2), \qquad (z_1,z_2), (\zeta_1,\zeta_2) \in \cc^{2n} \times \cc^{2n},$$
with $\sigma$ the canonical symplectic form on $\cc^{2n}$ defined in (\ref{11}). The set
\begin{equation}\label{graph1}
\widetilde{\lambda_{\mathcal{T}}}=\{(z_1,z_2,\zeta_1,\zeta_2) \in \cc^{4n} : (z_1,\zeta_1,z_2,\zeta_2) \in \lambda_{\mathcal{T}}\} \subset \cc^{4n},
\end{equation}
is then a non-negative Lagrangian plane of $\cc^{4n}$ equipped with the symplectic form
$$\sigma((z,\zeta),(\tilde{z},\tilde{\zeta}))=\langle \zeta, \tilde{z} \rangle -\langle z, \tilde{\zeta}\rangle=\sum_{j=1}^{2n}(\zeta_j \tilde{z}_j-z_j \tilde{\zeta}_j),$$
with $z=(z_1,...,z_{2n})$, $\tilde{z}=(\tilde{z}_1,....,\tilde{z}_{2n})$, $\zeta=(\zeta_1,...,\zeta_{2n})$, $\tilde{\zeta}=(\tilde{\zeta}_1,...,\tilde{\zeta}_{2n}) \in \cc^{2n}$. 
According to~\cite{mehler} (Proposition~5.1 and Proposition~5.5), there exists a complex-valued quadratic form 
\begin{equation}\label{pform}
p(x,y,\theta) = \langle (x,y, \theta), P (x,y, \theta) \rangle, \qquad (x,y) \in \rr^{2n}, \ \theta \in \rr^N, 
\end{equation}
where 
\begin{equation}\label{Pmatrix}
P=\left(
\begin{array}{ll}
P_{x,y;x,y} & P_{x,y; \theta} \\
P_{\theta; x,y} & P_{\theta; \theta} 
\end{array}
\right) \in \cc^{(2n+N) \times (2n+N)},
\end{equation}
is a symmetric matrix satisfying the conditions:

\begin{enumerate}
\item[$(i)$]  $\textrm{Im }P \geq 0$; 
\item[$(ii)$] The row vectors of  the submatrix
\begin{equation*}
\left(
\begin{array}{ll}
P_{\theta; x,y} & P_{\theta; \theta} 
\end{array}
\right) \in \cc^{N \times (2n+N)},
\end{equation*}
are linearly independent over $\cc$,
\end{enumerate}
parametrizing the non-negative Lagrangian plane
$$\widetilde{\lambda_{\mathcal{T}}}=\Big\{\Big(x,y,\frac{\partial p}{\partial x}(x,y,\theta),\frac{\partial p}{\partial y}(x,y,\theta)\Big) : \frac{\partial p}{\partial \theta}(x,y,\theta)=0\Big\}.$$
By using some integrations by parts as in~\cite[p.~442]{mehler} (see also Proposition~4.2 in~\cite{wahlberg}), this quadratic form $p$ allows to define the  tempered distribution
\begin{equation}\label{oscillint1}
K_{\mathcal{T}}(x,y) =\frac{1}{(2\pi)^{\frac{n+N}{2}}}\sqrt{\det\left(
\begin{array}{ll}
-i p''_{\theta,\theta} & p''_{\theta,y} \\
p''_{x,\theta} & ip''_{x,y} 
\end{array}
\right)} \int_{\rr^N} e^{i p(x,y,\theta)} \, d \theta \in \mathscr{S}'(\rr^{2n}), 
\end{equation}
as an oscillatory integral. Notice here that we do not prescribe the sign of the square root so the tempered distribution $K_{\mathcal{T}}$ is only determined up to its sign. Apart from this sign uncertainty, it is checked in~\cite[p.~444]{mehler} that this definition only depends on the non-negative complex symplectic linear transformation $\mathcal{T}$, and not on the choice of the parametrization of the non-negative Lagrangian plane $\widetilde{\lambda_{\mathcal{T}}}$ by the quadratic form $p$.

Associated to the non-negative complex symplectic linear transformation~$\mathcal{T}$ is therefore the Fourier integral operator 
$$\mathscr{K}_{\mathcal{T}} : \mathscr{S}(\rr^n) \rightarrow \mathscr{S}'(\rr^n),$$ 
defined by the kernel $K_{\mathcal{T}} \in \mathscr{S}'(\rr^{2n})$ as
$$\forall u,v \in \mathscr{S}(\rr^n), \quad \langle \mathscr{K}_{\mathcal{T}}u,v \rangle_{\mathscr{S}'(\rr^n),\mathscr{S}(\rr^n)}=\langle K_{\mathcal{T}},u \otimes v \rangle_{\mathscr{S}'(\rr^{2n}),\mathscr{S}(\rr^{2n})}.$$
It is proved in~\cite[p.~446]{mehler} that the adjoint operator 
$$\mathscr{K}_{\mathcal{T}}^* : \mathscr{S}(\rr^n) \rightarrow \mathscr{S}'(\rr^n),$$ 
defined as 
$$\forall u,v \in \mathscr{S}(\rr^n), \quad \langle \mathscr{K}_{\mathcal{T}}^*u,\overline{v} \rangle_{\mathscr{S}'(\rr^n),\mathscr{S}(\rr^n)}=\overline{\langle \mathscr{K}_{\mathcal{T}}v,\overline{u} \rangle}_{\mathscr{S}'(\rr^n),\mathscr{S}(\rr^n)},$$
is the Fourier integral operator $\mathscr{K}_{\overline{\mathcal{T}}^{-1}}$ associated to the non-negative complex symplectic linear transformation 
$$\overline{\mathcal{T}}^{-1} : \cc^{2n} \rightarrow \cc^{2n}.$$ 
Furthermore, the operator $\mathscr{K}_{\mathcal{T}}$ satisfies the Egorov formula proved in~\cite[p.~445]{mehler},
\begin{equation}\label{df1}
\forall u \in \mathscr{S}(\rr^n), \quad \big(\langle x_0,D_x \rangle-\langle \xi_0,x \rangle\big)\mathscr{K}_{\mathcal{T}}u=\mathscr{K}_{\mathcal{T}}\big(\langle y_0,D_x \rangle-\langle \eta_0,x \rangle\big)u,
\end{equation} 
with $(x_0,\xi_0)=\mathcal{T}(y_0,\eta_0)$. Thanks to this Egorov formula, it is proved in~\cite{mehler} (Proposition~5.8) that the operator $\mathscr{K}_{\mathcal{T}}$ is actually a continuous linear map on the Schwartz space $\mathscr{S}(\rr^n)$,
$$\mathscr{K}_{\mathcal{T}} : \mathscr{S}(\rr^n) \rightarrow \mathscr{S}(\rr^n).$$  
The mapping is then extended by duality for all $u \in \mathscr{S}'(\rr^n)$, $v \in \mathscr{S}(\rr^n)$,  
\begin{equation}\label{df1.2}
\langle \mathscr{K}_{\mathcal{T}}u,\overline{v} \rangle_{\mathscr{S}'(\rr^n),\mathscr{S}(\rr^n)}=\langle u,\overline{\mathscr{K}_{\mathcal{T}}^*v} \rangle_{\mathscr{S}'(\rr^n),\mathscr{S}(\rr^n)}=\langle u,\overline{\mathscr{K}_{\overline{\mathcal{T}}^{-1}}v} \rangle_{\mathscr{S}'(\rr^n),\mathscr{S}(\rr^n)},
\end{equation}
as a continuous linear map on the space of tempered distributions  $\mathscr{S}'(\rr^n)$,
$$\mathscr{K}_{\mathcal{T}} : \mathscr{S}'(\rr^n) \rightarrow \mathscr{S}'(\rr^n).$$
With this definition, the Egorov formula (\ref{df1}) extends by duality for tempered distributions
\begin{equation}\label{df1.1}
\forall u \in \mathscr{S}'(\rr^n), \quad \big(\langle x_0,D_x \rangle-\langle \xi_0,x \rangle\big)\mathscr{K}_{\mathcal{T}}u=\mathscr{K}_{\mathcal{T}}\big(\langle y_0,D_x \rangle-\langle \eta_0,x \rangle\big)u,
\end{equation} 
with $(x_0,\xi_0)=\mathcal{T}(y_0,\eta_0)$. Indeed, with $(\tilde{x}_0,\tilde{\xi}_0)=\overline{\mathcal{T}}^{-1}(\overline{x_0},\overline{\xi_0})$, we deduce from  (\ref{df1}) and (\ref{df1.2}) that for all $u \in \mathscr{S}'(\rr^n)$, $v \in \mathscr{S}(\rr^n)$,  
\begin{multline*}
\big\langle (\langle x_0,D_x \rangle-\langle \xi_0,x \rangle)\mathscr{K}_{\mathcal{T}}u,\overline{v} \big\rangle_{\mathscr{S}',\mathscr{S}}
=\big\langle \mathscr{K}_{\mathcal{T}}u,\overline{(\langle \overline{x_0},D_x \rangle-\langle \overline{\xi_0},x \rangle)v} \big\rangle_{\mathscr{S}',\mathscr{S}}\\
=\langle u,\overline{\mathscr{K}_{\overline{\mathcal{T}}^{-1}}(\langle \overline{x_0},D_x \rangle-\langle \overline{\xi_0},x \rangle)v} \rangle_{\mathscr{S}',\mathscr{S}}
=\langle u,\overline{(\langle \tilde{x}_0,D_x \rangle-\langle \tilde{\xi}_0,x \rangle)\mathscr{K}_{\overline{\mathcal{T}}^{-1}}v} \rangle_{\mathscr{S}',\mathscr{S}}\\
=\langle (\langle \overline{\tilde{x}_0},D_x \rangle-\langle \overline{\tilde{\xi}_0},x \rangle)u,\overline{\mathscr{K}_{\overline{\mathcal{T}}^{-1}}v} \rangle_{\mathscr{S}',\mathscr{S}}
=\big\langle \mathscr{K}_{\mathcal{T}}(\langle \overline{\tilde{x}_0},D_x \rangle-\langle \overline{\tilde{\xi}_0},x \rangle)u,\overline{v} \big\rangle_{\mathscr{S}',\mathscr{S}},
\end{multline*}
that is,
\begin{equation}\label{df1.3}
\forall  u \in \mathscr{S}'(\rr^n), \quad  \big(\langle x_0,D_x \rangle-\langle \xi_0,x \rangle\big)\mathscr{K}_{\mathcal{T}}u
= \mathscr{K}_{\mathcal{T}}\big(\langle y_0,D_x \rangle-\langle \eta_0,x \rangle\big)u,
\end{equation}
with $(x_0,\xi_0)=\mathcal{T}(y_0,\eta_0)$.
On the other hand, we recall from~\cite{mehler} that 
$$\mathscr{K}_{\mathcal{T}} : L^2(\rr^n) \rightarrow L^2(\rr^n),$$
defines a bounded operator on $L^2(\rr^n)$ whose operator norm satisfies
$$\|\mathscr{K}_{\mathcal{T}}\|_{\mathcal{L}(L^2(\rr^n))} \leq 1.$$
Indeed, it is proved in~\cite{mehler} (Proposition~5.12) that the operator $\mathscr{K}_{\mathcal{T}}$ is equal to a finite product of strongly continuous contraction semigroups on $L^2(\rr^n)$ at time $t=1$ generated by maximally dissipative quadratic operators $iQ_j^w(x,D_x)$,
$$\mathscr{K}_{\mathcal{T}}=e^{iQ_1^w(x,D_x)}...e^{iQ_k^w(x,D_x)},$$
where $Q_j$ are quadratic forms whose imaginary parts are non-negative $\textrm{Im }Q_j \geq 0$. It is also shown in~\cite{mehler} (Proposition~5.12) that the operator  $\mathscr{K}_{\mathcal{T}} : L^2(\rr^n) \rightarrow L^2(\rr^n),$ is invertible if and only if $\mathcal{T}$ is a real symplectic linear transformation. In this case, the operator  $\mathscr{K}_{\mathcal{T}}$ is a metaplectic operator associated to the real symplectic linear transformation $\mathcal{T}$ and the operator 
$$\mathscr{K}_{\mathcal{T}} : L^2(\rr^n) \rightarrow L^2(\rr^n),$$
defines a bijective isometry on $L^2(\rr^n)$.

The properties of this class of Fourier integral operators is summarized in the following proposition:

\bigskip

\begin{proposition}\label{FIOdef}
Associated to any non-negative complex symplectic linear transformation $\mathcal{T}$ is a Fourier integral operator 
$$\mathscr{K}_{\mathcal{T}} : \mathscr{S}(\rr^n) \rightarrow \mathscr{S}'(\rr^n),$$
whose kernel\footnote{determined up to its sign.} is the tempered distribution $K_{\mathcal{T}} \in \mathscr{S}'(\rr^{2n})$ defined in (\ref{oscillint1}), and whose adjoint 
$$\mathscr{K}_{\mathcal{T}}^*=\mathscr{K}_{\overline{\mathcal{T}}^{-1}} : \mathscr{S}(\rr^n) \rightarrow \mathscr{S}'(\rr^n),$$ 
is the Fourier integral operator associated to the non-negative complex symplectic linear transformation~$\overline{\mathcal{T}}^{-1}$. The Fourier integral operator $\mathscr{K}_{\mathcal{T}}$ defines a continuous mapping on the Schwartz space 
$$\mathscr{K}_{\mathcal{T}} : \mathscr{S}(\rr^n) \rightarrow \mathscr{S}(\rr^n),$$
which extends by duality as a continuous linear map on the space of tempered distributions 
$$\mathscr{K}_{\mathcal{T}} : \mathscr{S}'(\rr^n) \rightarrow \mathscr{S}'(\rr^n),$$
satisfying the Egorov formula
$$\forall (y_0,\eta_0) \in \cc^{2n}, \forall  u \in \mathscr{S}'(\rr^n), \quad  (\langle x_0,D_x \rangle-\langle \xi_0,x \rangle)\mathscr{K}_{\mathcal{T}}u
= \mathscr{K}_{\mathcal{T}}(\langle y_0,D_x \rangle-\langle \eta_0,x \rangle)u,$$
with $(x_0,\xi_0)=\mathcal{T}(y_0,\eta_0)$.
Furthermore, the Fourier integral operator 
$$\mathscr{K}_{\mathcal{T}} : L^2(\rr^n) \rightarrow L^2(\rr^n),$$  
is a bounded operator on $L^2(\rr^n)$ whose operator norm satisfies $\|\mathscr{K}_{\mathcal{T}}\|_{\mathcal{L}(L^2)} \leq 1$.
\end{proposition}

\bigskip

\begin{remark}
The kernel $K_{\mathcal{T}} \in \mathscr{S}'(\rr^{2n})$ of the Fourier integral operator $\mathscr{K}_{\mathcal{T}}$ appearing in the statement of  Proposition~\ref{FIOdef} is only determined up to its sign. In many cases as for the study of propagation of Gabor singularites in this work, this sign uncertainty is not an issue.  
\end{remark}

\section{Proofs of the main results}\label{proofth}

This section is devoted to the proofs of Theorems~\ref{th1},~\ref{th2} and~\ref{th3}. 
We begin by establishing the existence and uniqueness of evolution systems appearing in the statement of Theorem~\ref{th1}.

\subsection{Existence and uniqueness of evolution systems} 
Let $T>0$ and $q_t : \rr^{2n} \rightarrow \cc$ be a time-dependent complex-valued quadratic form 
$$q_t(x,\xi)=\sum_{\substack{\alpha, \beta \in \mathbb{N}^n \\ |\alpha+\beta| = 2}}(q_t)_{\alpha,\beta}x^{\alpha} \xi^{\beta},$$
with a non-positive real part 
\begin{equation}\label{sign}
\forall 0 \leq t \leq T, \quad \textrm{Re }q_t \leq 0,
\end{equation} 
and whose coefficients $(q_t)_{\alpha,\beta}$ depend continuously on the time variable $0 \leq t \leq T$.

This section is devoted to the proof of existence and uniqueness of an evolution system for the non-autonomous Cauchy problem
$$\left\lbrace\begin{array}{l}
\frac{du(t)}{dt}=q_t^w(x,D_x)u(t), \qquad 0 \leq \tau <t \leq T,\\
u(\tau)=v.
\end{array}\right.
$$
We follow the theory of non-autonomous evolution systems developed in~\cite{pazy} (Chapter~5).

According to~\cite[pp. 425-426]{mehler}, the assumption (\ref{sign}) implies that $(q_t^w(x,D_x))_{0 \leq t \leq T}$ is a family of infinitesimal generators of strongly continuous contraction semigroups on $L^2(\rr^n)$. This family $(q_t^w(x,D_x))_{0 \leq t \leq T}$ is therefore stable~\cite[p.~131]{pazy} in the sense of Definition~5.2.1 in~\cite{pazy}. Let $B$ be the Hilbert space defined in (\ref{hilbert}). The space $B$ contains the Schwartz space $\mathscr{S}(\rr^n)$. This Hilbert space is therefore densely and continuously imbedded in $L^2(\rr^n)$,  
$$\forall u \in B, \quad \|u\|_{L^2(\rr^n)} \leq \|u\|_{B}.$$
It follows from (\ref{dom1}) that 
$$\forall t \geq 0, \quad B \subset D(q_t^w).$$
We observe that  the quadratic operator 
\begin{equation}\label{operator}
q_t^w(x,D_x)=\sum_{\substack{\alpha, \beta \in \mathbb{N}^n \\ |\alpha+\beta| = 2}}(q_t)_{\alpha,\beta}\frac{x^{\alpha}D_x^{\beta}+D_x^{\beta} x^{\alpha}}{2},
\end{equation}
satisfies for all $u \in B$,
\begin{multline*}
\|q_t^w(x,D_x)u\|_{L^2(\rr^n)} \leq \sum_{\substack{\alpha, \beta \in \mathbb{N}^n \\ |\alpha+\beta| = 2}}|(q_t)_{\alpha,\beta}|\Big(\|x^{\alpha}D_x^{\beta}u\|_{L^2(\rr^n)}+\frac{1}{2}
\|[D_x^{\beta},x^{\alpha}]u\|_{L^2(\rr^n)}\Big)\\
\leq \frac{3}{2}\Big(\sum_{\substack{\alpha, \beta \in \mathbb{N}^n \\ |\alpha+\beta| = 2}}|(q_t)_{\alpha,\beta}|\Big)\|u\|_{B}.
\end{multline*}
This implies that $q_t^w(x,D_x)$ defines a bounded operator from $B$ to $L^2(\rr^n)$, 
$$\|q_t^w(x,D_x)\|_{\mathcal{L}(B,L^2)} \leq \frac{3}{2}\sum_{\substack{\alpha, \beta \in \mathbb{N}^n \\ |\alpha+\beta| = 2}}|(q_t)_{\alpha,\beta}|,$$
so that the mapping 
$$t \in [0,T] \mapsto q_t^w(x,D_x) \in \big(\mathcal{L}(B,L^2), \|\cdot\|_{\mathcal{L}(B,L^2)}\big),$$ 
is continuous.

We now check that $B$ is $q_t^w$-admissible for all $0 \leq t \leq T$. We recall from~\cite[p.~122]{pazy} (Definition~4.5.3) that while denoting $(e^{\tau q_t^w})_{\tau \geq 0}$ the strongly continuous contraction semigroup generated by the quadratic operator $q_t^w(x,D_x)$, it means that 
\begin{equation}\label{lh2}
\forall  0 \leq t \leq T, \forall \tau \geq 0, \quad e^{\tau q_t^w}(B) \subset B
\end{equation} 
and that for all $0 \leq t \leq T$, the restriction of $(e^{\tau q_t^w})_{\tau \geq 0}$ to $B$ is a strongly continuous semigroup in $B$, that is, strongly continuous in the $B$-norm.
Let $0 \leq t \leq T$. We know from~\cite{mehler} (Theorem~5.12) that the strongly continuous contraction semigroup $e^{\tau q_t^w}$ at time $\tau \geq 0$ is equal to the Fourier integral operator
\begin{equation}\label{lh0}
e^{\tau q_t^w} =\mathscr{K}_{e^{2i \tau F_t}},
\end{equation} 
associated to the non-negative complex symplectic linear transformation 
$$e^{2i \tau F_t} : \cc^{2n} \rightarrow \cc^{2n}.$$ 
We deduce from Proposition~\ref{FIOdef} and (\ref{lh0}) that for all $(x_1,\xi_1) \in \rr^{2n}$, $(x_2,\xi_2) \in \rr^{2n}$, $0 \leq t \leq T$, $\tau \geq 0$, $u \in \mathscr{S}'(\rr^n)$,
\begin{equation}\label{lh1.1}
\langle (-\xi_1,x_1),(x,D_x)\rangle  e^{\tau q_t^w}u= e^{\tau q_t^w}\langle (-\sigma)e^{-2i\tau F_t}(x_1,\xi_1),(x,D_x)\rangle u
\end{equation}
and
\begin{multline}\label{lh1}
\langle (-\xi_1,x_1),(x,D_x)\rangle  \langle (-\xi_2,x_2),(x,D_x) \rangle e^{\tau q_t^w}u\\
= e^{\tau q_t^w}\langle \sigma e^{-2i\tau F_t}(x_1,\xi_1),(x,D_x)\rangle  \langle \sigma e^{-2i \tau F_t}(x_2,\xi_2),(x,D_x) \rangle u,
\end{multline}
with $\sigma=\left(\begin{array}{cc}0 & I_n \\ -I_n & 0\end{array}\right)$. 
With $\|\cdot\|$ the Euclidean norm on $\cc^n$, we notice that 
\begin{equation}\label{lh-1}
\|\langle (a,b),(x,D_x) \rangle u\|_{L^2}
 \leq \|(a,b)\|\sum_{j=1}^n\big(\|x_ju\|_{L^2}+\|D_{x_j}u\|_{L^2}\big)\leq 2n\|(a,b)\|\|u\|_B.
\end{equation}
On the other hand, we deduce from the estimates (\ref{lh-1}) that 
\begin{align}\label{lh-2}
& \ \|\langle (a_1,b_1),(x,D_x) \rangle \langle (a_2,b_2),(x,D_x) \rangle u\|_{L^2} \\  \notag
\leq & \ \|(a_1,b_1)\|\sum_{j=1}^n\big(\|x_j \langle (a_2,b_2),(x,D_x) \rangle u\|_{L^2}+\|D_{x_j} \langle (a_2,b_2),(x,D_x) \rangle u\|_{L^2}\big)\\ \notag
\leq & \ \|(a_1,b_1)\|\\ \notag
    \times & \sum_{j=1}^n\big(\|\langle (a_2,b_2),(x,D_x) \rangle x_ju\|_{L^2}+2\|(a_2,b_2)\|\|u\|_{L^2}+\| \langle (a_2,b_2),(x,D_x) \rangle D_{x_j}u\|_{L^2}\big)\\ \notag
\leq & \ \|(a_1,b_1)\|\|(a_2,b_2)\|\\ \notag
    \times & \Big(\sum_{1 \leq j,k \leq n}\big(\|x_k x_ju\|_{L^2}+\|D_{x_k} x_ju\|_{L^2}+\|D^2_{x_k,x_j}u\|_{L^2}+\| x_k D_{x_j}u\|_{L^2})+2n\|u\|_{L^2}\Big).
\end{align}
We obtain from (\ref{lh-2}) that there exists a positive constant $C_n>0$ such that 
\begin{align}\label{lh-3}
&\ \|\langle (a_1,b_1),(x,D_x) \rangle \langle (a_2,b_2),(x,D_x) \rangle u\|_{L^2} \\ \notag
\leq & \ \|(a_1,b_1)\|\|(a_2,b_2)\|\\ \notag
& \ \qquad \times \Big(\sum_{1 \leq j,k \leq n}\big(\|x_k x_ju\|_{L^2}+\|D^2_{x_k,x_j}u\|_{L^2}+2\| x_k D_{x_j}u\|_{L^2})+3n\|u\|_{L^2}\Big)\\ \notag
\leq & \ (4n+3)n\|(a_1,b_1)\|\|(a_2,b_2)\|\|u\|_B.
\end{align}
It follows from (\ref{lh1.1}), (\ref{lh1}), (\ref{lh-1}) and (\ref{lh-3}) that the inclusion (\ref{lh2}) holds and there exists a positive constant $C>0$ such that 
\begin{equation}\label{lh3}
\forall 0 \leq t \leq T, \forall \tau \geq 0, \forall u \in B, \quad \|e^{\tau q_t^w}u\|_{B}  \leq C e^{4\tau \|F_t\|}\|u\|_B,
\end{equation} 
since $(e^{\tau q_t^w})_{\tau  \geq 0}$ is a strongly continuous contraction semigroup on $L^2(\rr^n)$. The operator $e^{\tau q_t^w}$ is therefore bounded on $B$ for all $0 \leq t \leq T$, $\tau \geq 0$. It remains to check that for all $0 \leq t \leq T$ and $u \in B$, the mapping $\tau \in [0,+\infty[ \mapsto e^{\tau q_t^w}u \in B$ is continuous. It is sufficient to prove that for all $0 \leq t \leq T$, $(x_1,\xi_1) \in \rr^{2n}$, $(x_2,\xi_2) \in \rr^{2n}$ and $u \in B$, the mappings 
\begin{equation}\label{lh6}
\tau \in [0,+\infty[ \mapsto \langle (-\xi_1,x_1),(x,D_x)\rangle  e^{\tau q_t^w}u \in L^2(\rr^n)
\end{equation} 
and 
\begin{equation}\label{lh7}
\tau  \in [0,+\infty[ \mapsto \langle (-\xi_1,x_1),(x,D_x)\rangle  \langle (-\xi_2,x_2),(x,D_x) \rangle e^{\tau q_t^w}u \in L^2(\rr^n),
\end{equation}
are continuous. For all $\tau, \tau_0 \geq 0$, we deduce from (\ref{lh1.1}) that 
\begin{align}\label{lh4}
& \ \|\langle (-\xi_1,x_1),(x,D_x)\rangle  e^{\tau q_t^w}u-\langle (-\xi_1,x_1),(x,D_x)\rangle  e^{\tau_0q_t^w}u\|_{L^2}\\ \notag
= & \ \|e^{\tau q_t^w}\langle \sigma e^{-2i\tau F_t}(x_1,\xi_1),(x,D_x)\rangle u-e^{\tau_0q_t^w}\langle \sigma e^{-2i\tau_0F_t}(x_1,\xi_1),(x,D_x)\rangle u\|_{L^2}\\ \notag
\leq & \  \|e^{\tau q_t^w}\langle \sigma (e^{-2i\tau F_t}-e^{-2i\tau_0F_t})(x_1,\xi_1),(x,D_x)\rangle u\|_{L^2}\\ \notag
& \ +\|(e^{\tau q_t^w}-e^{\tau_0q_t^w})\langle \sigma e^{-2i\tau_0F_t}(x_1,\xi_1),(x,D_x)\rangle u\|_{L^2}.
\end{align}
By using that $(e^{\tau q_t^w})_{\tau \geq 0}$ is a strongly continuous contraction semigroup on $L^2(\rr^n)$, we obtain from (\ref{lh-1}) that there exists a positive constant $C>0$ such that for all $u \in B$, $\tau,\tau_0 \geq 0$, $0 \leq t \leq T$,
\begin{align}\label{lh5}
& \ \|\langle (-\xi_1,x_1),(x,D_x)\rangle  e^{\tau q_t^w}u-\langle (-\xi_1,x_1),(x,D_x)\rangle  e^{\tau_0q_t^w}u\|_{L^2}\\ \notag
\leq & \  \|\langle \sigma (e^{-2i\tau F_t}-e^{-2i\tau_0F_t})(x_1,\xi_1),(x,D_x)\rangle u\|_{L^2} \\ \notag
& \ +\|(e^{\tau q_t^w}-e^{\tau_0q_t^w})\langle \sigma e^{-2i\tau_0F_t}(x_1,\xi_1),(x,D_x)\rangle u\|_{L^2}\\ \notag
\leq & \  C\|e^{-2i\tau F_t}-e^{-2i\tau_0F_t}\| \|u\|_{B}  +\|(e^{\tau q_t^w}-e^{\tau_0q_t^w})\underbrace{\langle \sigma e^{-2i\tau_0F_t}(x_1,\xi_1),(x,D_x)\rangle u}_{\in L^2(\rr^n)}\|_{L^2}.
\end{align}
Then, the continuity of the mapping (\ref{lh6}) follows from the continuity of the mapping $\tau \in [0,+\infty[ \mapsto  e^{\tau q_t^w}v \in L^2(\rr^n)$ for $v \in L^2(\rr^n)$. The very same arguments allow to prove the continuity of the mapping (\ref{lh7}). It proves that $B$ is $q_t^w$-admissible for all $0 \leq t \leq T$. It follows from~\cite{pazy} (Definition~1.10.3 and Theorem~4.5.5) that the part of the operator $q_t^w(x,D_x)$ in $B$, that is, the operator
$$\begin{array}{cl}
\widetilde{q_t^w(x,D_x)} :  \{ u \in B \cap D(q_t^w) : q_t^wu \in B\} & \rightarrow B\\
 u & \mapsto q_t^w(x,D_x)u,\end{array}$$ 
is the infinitesimal generator of a strongly continuous semigroup on $B$. Furthermore, this strongly continuous semigroup on $B$ is given by the restriction of $L^2$-semigroup $(e^{\tau q_t^w})_{\tau \geq 0}$ to $B$,
\begin{equation}\label{lh6.99}
\forall 0 \leq t \leq T, \forall \tau \geq 0, \forall u \in B, \quad e^{\tau \widetilde{q_t^w}}u=e^{\tau q_t^w}u.
\end{equation}
We deduce from (\ref{lh3}) that the strongly continuous semigroup $(e^{\tau \tilde{q_t^w}})_{\tau \geq 0}$ on $B$ satisfies   
\begin{equation}\label{lh3.1}
\forall 0 \leq t \leq T,  \forall \tau \geq 0, \quad \|e^{\tau \widetilde{q_t^w}}\|_{\mathcal{L}(B)}  \leq C e^{4\tau \|F_t\|}.
\end{equation} 
It follows from~\cite{pazy} (Theorem~1.5.3) that the resolvent set of the operator $\widetilde{q_t^w(x,D_x)}$ contains the ray 
\begin{equation}\label{ray}
]4\|F_t\|,+\infty[.
\end{equation}
Recalling the continuity of the mapping $t \in [0,T] \mapsto F_t=\sigma Q_t \in M_{2n}(\cc)$,
we set
\begin{equation}\label{lh7.99}
0 \leq \omega=\sup_{0 \leq t \leq T}\|F_t\| <+\infty.
\end{equation}
Let $k \geq 1$ and $0 \leq t_1 \leq t_2 \leq ... \leq t_k \leq T$ and $\tau_1, ...., \tau_k \geq 0$.   
We deduce from (\ref{lh1.1}) and (\ref{lh6.99}) that for all $(x_1,\xi_1) \in \rr^{2n}$, $(x_2,\xi_2) \in \rr^{2n}$, $0 \leq t \leq T$, $\tau \geq 0$, $u \in B$,
\begin{multline}\label{lh1.12}
\langle (x_1,\xi_1),(x,D_x)\rangle  e^{\tau_1\widetilde{q_{t_1}^w}}...e^{\tau_k\widetilde{q_{t_k}^w}}u\\
=(-1)^ke^{\tau_1\widetilde{q_{t_1}^w}}...e^{\tau_k\widetilde{q_{t_k}^w}}\langle \sigma e^{-2i\tau_kF_{t_k}}\sigma...\sigma e^{-2i\tau_1F_{t_1}}\sigma(x_1,\xi_1),(x,D_x)\rangle u
\end{multline}
and
\begin{align}\label{lh1.22}
& \ \langle (x_1,\xi_1),(x,D_x)\rangle \langle (x_2,\xi_2),(x,D_x)\rangle e^{\tau_1\widetilde{q_{t_1}^w}}...e^{\tau_k\widetilde{q_{t_k}^w}}u\\ \notag
= & \ e^{\tau_1\widetilde{q_{t_1}^w}}...e^{\tau_k\widetilde{q_{t_k}^w}}\langle \sigma e^{-2i\tau_kF_{t_k}}\sigma...\sigma e^{-2i\tau_1F_{t_1}}\sigma(x_1,\xi_1),(x,D_x)\rangle \\ \notag
& \ \qquad \qquad \qquad \qquad \qquad \qquad \qquad \quad  \langle \sigma e^{-2i\tau_kF_{t_k}}\sigma...\sigma e^{-2i\tau_1F_{t_1}}\sigma(x_2,\xi_2),(x,D_x) \rangle u,
\end{align}
with 
$$\sigma=\left(\begin{array}{cc}0 & I_n \\ -I_n & 0\end{array}\right).$$
We observe from (\ref{lh7.99}) that 
\begin{multline}\label{cc1}
\|\sigma e^{-2i\tau_kF_{t_k}}\sigma...\sigma e^{-2i\tau_1F_{t_1}}\sigma(x_j,\xi_j)\| \leq e^{2(\tau_k\|F_{t_k}\|+...+\tau_1\|F_{t_1}\|)}\|(x_j,\xi_j)\|\\
\leq e^{2(\tau_1+...+\tau_k)\omega}\|(x_j,\xi_j)\|,
\end{multline}
since $\|\sigma\|=1$.
Recalling that $\|e^{\tau_jq_{t_j}^w}\|_{\mathcal{L}(L^2)} \leq 1$, we deduce from (\ref{lh-1}), (\ref{lh-3}), (\ref{lh1.12}), (\ref{lh1.22}) and (\ref{cc1}) that 
\begin{multline}\label{lh1.123}
\|\langle (x_1,\xi_1),(x,D_x)\rangle  e^{\tau_1\widetilde{q_{t_1}^w}}...e^{\tau_k\widetilde{q_{t_k}^w}}u\|_{L^2} \\ \leq \|\langle \sigma e^{-2i\tau_kF_{t_k}}\sigma...\sigma e^{-2i\tau_1F_{t_1}}\sigma(x_1,\xi_1),(x,D_x)\rangle u\|_{L^2} \leq 2n e^{2(\tau_1+...+\tau_k)\omega}\|(x_1,\xi_1)\|\|u\|_B
\end{multline}
and
\begin{align}\label{lh1.223}
& \ \|\langle (x_1,\xi_1),(x,D_x)\rangle \langle (x_2,\xi_2),(x,D_x)\rangle e^{\tau_1\widetilde{q_{t_1}^w}}...e^{\tau_k\widetilde{q_{t_k}^w}}u\|_{L^2}\\ \notag
\leq &\   \|\langle \sigma e^{-2i\tau_kF_{t_k}}\sigma...\sigma e^{-2i\tau_1F_{t_1}}\sigma(x_1,\xi_1),(x,D_x)\rangle \\ \notag
& \ \qquad \qquad \qquad \qquad \qquad \qquad \qquad  \langle \sigma e^{-2i\tau_kF_{t_k}}\sigma...\sigma e^{-2i\tau_1F_{t_1}}\sigma(x_2,\xi_2),(x,D_x) \rangle u\|_{L^2}\\ \notag
\leq & \ (4n+3)n e^{4(\tau_1+...+\tau_k)\omega}\|(x_1,\xi_1)\|\|(x_2,\xi_2)\|\|u\|_B.
\end{align}
We deduce from (\ref{lh1.123}) and (\ref{lh1.223}) that there exists a positive constant $M \geq 1$ such that for all $k \geq 1$, $0 \leq t_1 \leq t_2 \leq ... \leq t_k \leq T$ and $\tau_1, ...., \tau_k \geq 0$,
$$\|e^{\tau_1\widetilde{q_{t_1}^w}}...e^{\tau_k\widetilde{q_{t_k}^w}}\|_{\mathcal{L}(B)} \leq Me^{4(\tau_1+...+\tau_k)\omega}.$$
According to (\ref{ray}) and (\ref{lh7.99}), it follows from~\cite[p.~131]{pazy} (Theorem~5.2.2) that the family of generators $(\widetilde{q_{t}^w})_{0 \leq t \leq T}$ is stable in $B$. The family of operators $(q_{t}^w(x,D_x))_{0 \leq t \leq T}$ satisfies the assumptions of Theorem~5.3.1 in~\cite[p.~135]{pazy}. We deduce from this result that there exists a unique evolution system $(\mathscr{U}(t,\tau))_{0 \leq \tau \leq t \leq T}$ in $L^2(\rr^n)$ satisfying
\begin{equation}\label{cc2}
\forall 0 \leq \tau \leq t \leq T, \quad \|\mathscr{U}(t,\tau)\| \leq 1,
\end{equation} 
\begin{equation}\label{cc3}
\forall 0 \leq \tau \leq T, \forall v \in B, \quad \frac{\partial^+}{\partial t}\mathscr{U}(t,\tau)v|_{t=\tau}=q_\tau^w(x,D_x)v,
\end{equation} 
\begin{equation}\label{cc4}
\forall 0 \leq \tau \leq t \leq T, \forall v \in B, \quad \frac{\partial}{\partial \tau}\mathscr{U}(t,\tau)v=-\mathscr{U}(t,\tau)q_\tau^w(x,D_x)v,
\end{equation} 
where the derivative from the right in (\ref{cc3}) and the derivative in (\ref{cc4}) are in the strong sense in $L^2(\rr^n)$.

\subsection{Existence and uniqueness of $B$-valued solutions}
We consider the notion of $B$-valued solutions given in Definition~\ref{def1}.
The existence of the evolution system given in the previous section is actually not sufficient to prove the existence of $B$-valued solutions to the non-autonomous Cauchy problem
\begin{equation}\label{cauchy}
\left\lbrace\begin{array}{l}
\frac{du(t)}{dt}=q_t^w(x,D_x)u(t), \qquad 0 \leq \tau <t \leq T,\\
u(\tau)=v.
\end{array}\right.
\end{equation}
However, we already know from~\cite{pazy} (Theorem~5.4.2) that if the non-autonomous Cauchy problem (\ref{cauchy}) has a $B$-valued solution $u$ then this solution is unique and given by the following formula
\begin{equation}\label{kk2}
u(t)=\mathscr{U}(t,\tau)v, \qquad 0 \leq \tau \leq t \leq T.
\end{equation}
Indeed, the existence of the evolution system $(\mathscr{U}(t,\tau))_{0 \leq \tau \leq t \leq T}$ only ensures the uniqueness of $B$-valued solutions 
but not the existence of $B$-valued solutions as the function 
$$u(t)=\mathscr{U}(t,\tau)v,$$ 
is not in general a $B$-valued solution. In fact, the subspace $B$ does not need to be an invariant subspace for $\mathscr{U}(t,\tau)$, and even if it is such an invariant subspace, the mapping $t \mapsto \mathscr{U}(t,\tau)v$ for $v \in B$ does not need to be continuous in the $B$-norm.

We now study the existence of $B$-valued solutions for the non-autonomous Cauchy problem (\ref{cauchy}). Setting 
$$\mathscr{H}=-\Delta_x+x^2,$$
this harmonic oscillator defines an isomorphism from $B$ onto $L^2(\rr^n)$. Furthermore, we observe that its Weyl symbol belongs to the following Shubin class 
$$\xi^2+x^2 \in S(\langle (x,\xi)\rangle^2,\langle (x,\xi)\rangle^{-2}(dx^2+d\xi^2)).$$ 
The H\"ormander notation 
$$S(\langle (x,\xi)\rangle^m,\langle (x,\xi)\rangle^{-2}(dx^2+d\xi^2)), \quad m \in \rr,$$
refers to the class of smooth complex-valued symbols satisfying the estimates
$$\forall \alpha, \beta \in \mathbb{N}^n, \exists C_{\alpha,\beta}>0, \forall (x,\xi) \in \rr^{2n}, \quad |\partial_x^{\alpha}\partial_{\xi}^{\beta}a(x,\xi)| \leq C_{\alpha,\beta}\langle (x,\xi)\rangle^{m-|\alpha|-|\beta|}.$$
We recall from~\cite{shubin} (Theorem~25.4) (see also~\cite{cappiello}) that the inverse of the harmonic oscillator $\mathscr{H}^{-1}$ writes as a pseudodifferential operator with a symbol belonging to the Shubin class 
$$S(\langle (x,\xi)\rangle^{-2},\langle (x,\xi)\rangle^{-2}(dx^2+d\xi^2)).$$
On the other hand, we notice from (\ref{operator}) that 
\begin{multline}\label{kk5}
\mathscr{H}q_t^w(x,D_x)\mathscr{H}^{-1}=\sum_{\substack{\alpha, \beta \in \mathbb{N}^n \\ |\alpha+\beta| = 2}}(q_t)_{\alpha,\beta}\mathscr{H}\frac{x^{\alpha}D_x^{\beta}+D_x^{\beta} x^{\alpha}}{2}\mathscr{H}^{-1}\\
=q_t^w(x,D_x)+\sum_{\substack{\alpha, \beta \in \mathbb{N}^n \\ |\alpha+\beta| = 2}}(q_t)_{\alpha,\beta}\Big[\mathscr{H},\frac{x^{\alpha}D_x^{\beta}+D_x^{\beta} x^{\alpha}}{2}\Big]\mathscr{H}^{-1}\\
=q_t^w(x,D_x)+\frac{1}{i}\sum_{\substack{\alpha, \beta \in \mathbb{N}^n \\ |\alpha+\beta| = 2}}(q_t)_{\alpha,\beta}\textrm{Op}^w\big(\{\xi^2+x^2,x^{\alpha}\xi^{\beta}\}\big)\mathscr{H}^{-1},
\end{multline}
where $\textrm{Op}^w\big(\{\xi^2+x^2,x^{\alpha}\xi^{\beta}\}\big)$ denotes the Weyl quantization of the Poisson bracket 
$$\big\{\xi^2+x^2,x^{\alpha}\xi^{\beta}\big\}=\sum_{j=1}^n\Big(\frac{\partial}{\partial\xi_j}(\xi^2+x^2)\frac{\partial}{\partial x_j}(x^{\alpha}\xi^{\beta})-\frac{\partial}{\partial x_j}(\xi^2+x^2)\frac{\partial}{\partial \xi_j}(x^{\alpha}\xi^{\beta})\Big).$$
We observe that this symbol belongs to the Shubin class 
$$S(\langle (x,\xi)\rangle^{2},\langle (x,\xi)\rangle^{-2}(dx^2+d\xi^2)).$$ 
By composition, we obtain that the Weyl symbol of the time-independent operator 
$$\textrm{Op}^w\big(\{\xi^2+x^2,x^{\alpha}\xi^{\beta}\}\big)\mathscr{H}^{-1},$$ 
belongs to the Shubin class
$$S(1,\langle (x,\xi)\rangle^{-2}(dx^2+d\xi^2)).$$ 
We therefore deduce from the Calder\'on-Vaillancourt theorem that 
$$t \in [0,T] \mapsto \sum_{\substack{\alpha, \beta \in \mathbb{N}^n \\ |\alpha+\beta| = 2}}(q_t)_{\alpha,\beta}\textrm{Op}^w\big(\{\xi^2+x^2,x^{\alpha}\xi^{\beta}\}\big)\mathscr{H}^{-1},$$
is a $L^2$-norm continuous (and thus also strongly continuous) family of bounded operator on $L^2(\rr^n)$. We can therefore apply~\cite{pazy} (Theorem~5.4.6) to obtain that the unique evolution system $(\mathscr{U}(t,\tau))_{0 \leq \tau \leq t \leq T}$ on $L^2(\rr^n)$ satisfying (\ref{cc2}), (\ref{cc3}) and (\ref{cc4}) also verifies
\begin{equation}\label{cc10}
\forall 0 \leq \tau \leq t \leq T, \quad \mathscr{U}(t,\tau)(B) \subset B
\end{equation}
and for all $v \in B$, the mapping $\mathscr{U}(t,\tau)v$ is continuous in $B$ for $0 \leq \tau \leq t \leq T$. We finally deduce from~\cite{pazy} (Theorem~5.4.3) that for all $v \in B$, $\mathscr{U}(t,\tau)v$ is the unique $B$-valued solution of the non-autonomous Cauchy problem (\ref{cauchy}). This ends the proof of Theorem~\ref{th1}.

\subsection{Some computations in the Weyl quantization}
This section is devoted to derive a formula for the Weyl symbol of the evolution operators. We begin with some symbolic computations in the Weyl quantization.

Let $T>0$ and $q_t : \rr^{2n} \rightarrow \cc$ be a time-dependent complex-valued quadratic form 
\begin{equation}\label{carre2}
q_t(x,\xi)=\sum_{\substack{\alpha, \beta \in \mathbb{N}^n \\ |\alpha+\beta| = 2}}(q_t)_{\alpha,\beta}x^{\alpha} \xi^{\beta},
\end{equation}
with a non-positive real part 
$$\forall 0 \leq t \leq T, \quad \textrm{Re }q_t \leq 0,$$ 
and whose coefficients $(q_t)_{\alpha,\beta}$ depend continuously on the time variable $0 \leq t \leq T$.
Let $Q_t \in \cc^{2n \times 2n}$ be the symmetric matrix defining the time-dependent quadratic form 
$$q_t(X)=\langle Q_t X,X \rangle, \quad 0 \leq t \leq T, \ X=(x,\xi) \in \rr^{2n}.$$ 
By assumption,
$\textrm{Re }Q_t \leq 0$ is a negative semidefinite symmetric matrix and the mapping $t \in [0,T] \mapsto Q_t \in \cc^{2n \times 2n}$ is a $C^0$ function on $[0,T]$.
Our ansatz is to find out a function
\begin{equation}\label{jk2} 
g_{t,\tau}(X)=\langle G_{t,\tau} X,X\rangle+h(t,\tau), \quad X=(x,\xi) \in \rr^{2n},
\end{equation}
with $G_{t,\tau} \in \cc^{2n \times 2n}$ a symmetric matrix depending continuously differentiably on $(t, \tau) \in [0,T]^2$ and $h(t,\tau)$ a continuously differentiable complex-valued function, satisfying the equations
\begin{equation}\label{jk3}
\frac{d}{dt}\big(e^{g_{t,\tau}}\big)=q_t \#^we^{g_{t,\tau}}, \qquad \frac{d}{d\tau}\big(e^{g_{t,\tau}}\big)=-e^{g_{t,\tau}}\#^wq_{\tau} ,
\end{equation}
where $a \#^w b$ denotes the Moyal product, that is, the symbol obtained by composition in the Weyl quantization
\begin{equation}\label{jk4}
(a \#^w b)(x,\xi)=\Big[e^{\frac{i}{2}\sigma(D_x,D_{\xi};D_y,D_{\eta})}\big(a(x,\xi)b(y,\eta)\big)\Big]\Big|_{(x,\xi)=(y,\eta)}.
\end{equation}
By using that $q_t$ is a quadratic symbol, we deduce from (\ref{jk3}) and (\ref{jk4}) that 
\begin{multline}\label{jk5}
\frac{\partial g_{t,\tau}}{\partial t}(X) e^{g_{t,\tau}(X)}=\Big[q_t(X)e^{g_{t,\tau}(Y)}+\frac{i}{2}\sigma(D_X;D_Y)\big(q_t(X)e^{g_{t,\tau}(Y)}\big)\\ +\frac{1}{2!}\Big(\frac{i}{2}\Big)^2\sigma(D_X;D_Y)^2\big(q_t(X)e^{g_{t,\tau}(Y)}\big)\Big]\Big|_{X=Y}
\end{multline} 
and
\begin{multline}\label{jk5.k}
\frac{\partial g_{t,\tau}}{\partial \tau}(X) e^{g_{t,\tau}(X)}=-\Big[e^{g_{t,\tau}(X)}q_{\tau}(Y)+\frac{i}{2}\sigma(D_X;D_Y)\big(e^{g_{t,\tau}(X)}q_{\tau}(Y)\big)\\ +\frac{1}{2!}\Big(\frac{i}{2}\Big)^2\sigma(D_X;D_Y)^2\big(e^{g_{t,\tau}(X)}q_{\tau}(Y)\big)\Big]\Big|_{X=Y},
\end{multline} 
with $X=(x,\xi) \in \rr^{2n}$ and $Y=(y,\eta) \in \rr^{2n}$.
Some direct computations provide that 
\begin{align}\label{jk6}
& \ \sigma(D_X;D_Y)\big(q_t(X)e^{g_{t,\tau}(Y)}\big)=-\langle \sigma \nabla_X, \nabla_Y\rangle \big(q_t(X)e^{g_{t,\tau}(Y)}\big)\\  \notag
= & \ -\langle \sigma \nabla_X q_t(X), \nabla_Yg_{t,\tau}(Y)\rangle e^{g_{t,\tau}(Y)}=-4\langle \sigma Q_tX, G_{t,\tau} Y\rangle e^{g_{t,\tau}(Y)}\\ \notag
=& \ -4\langle G_{t,\tau}\sigma Q_tX, Y\rangle e^{g_{t,\tau}(Y)}
\end{align}
and
\begin{align}\label{jk7}
& \ \sigma(D_X;D_Y)^2\big(q_t(X)e^{g_{t,\tau}(Y)}\big)\\ \notag
= & \ 4\sum_{1 \leq j,k \leq 2n}(\sigma \nabla_X)_j(\nabla_Y)_j\Big((G_{t,\tau} \sigma Q_tX)_k Y_k e^{g_{t,\tau}(Y)}\Big)\\ \notag
= & \ 4\sum_{1 \leq j \leq 2n}(\sigma \nabla_X)_j\Big((G_{t,\tau} \sigma Q_tX)_j\Big)e^{g_{t,\tau}(Y)}\\ \notag
& \ +8\sum_{1 \leq j,k \leq 2n}(\sigma \nabla_X)_j\Big((G_{t,\tau} \sigma Q_tX)_k\Big) Y_k(G_{t,\tau}Y)_j e^{g_{t,\tau}(Y)}.
\end{align}
While separating terms by homogeneity degree, we obtain from (\ref{jk5}), (\ref{jk5.k}), (\ref{jk6}) and (\ref{jk7}) the following equations
\begin{multline}\label{jk8}
\langle \partial_tG_{t,\tau} X,X\rangle=q_t(X)-2i\langle G_{t,\tau}\sigma Q_tX,X \rangle\\
-\sum_{1 \leq j,k \leq 2n}(\sigma \nabla_X)_j\Big((G_{t,\tau} \sigma Q_tX)_k\Big) X_k(G_{t,\tau}X)_j,
\end{multline}
\begin{multline}\label{jk8.k}
\langle \partial_{\tau}G_{t,\tau} X,X\rangle=-q_{\tau}(X)-2i\langle G_{t,\tau}\sigma Q_{\tau}X,X \rangle\\
+\sum_{1 \leq j,k \leq 2n}(\sigma \nabla_X)_j\Big((G_{t,\tau} \sigma Q_{\tau}X)_k\Big) X_k(G_{t,\tau}X)_j,
\end{multline}
\begin{equation}\label{jk9}
\partial_th(t,\tau)=-\frac{1}{2}\sum_{1 \leq j \leq 2n}(\sigma \nabla_X)_j\Big((G_{t,\tau} \sigma Q_tX)_j\Big),
\end{equation}
\begin{equation}\label{jk9.k}
\partial_{\tau}h(t,\tau)=\frac{1}{2}\sum_{1 \leq j \leq 2n}(\sigma \nabla_X)_j\Big((G_{t,\tau} \sigma Q_{\tau}X)_j\Big).
\end{equation}
We notice that 
\begin{align*}
& \ \sum_{1 \leq j,k \leq 2n}(\sigma \nabla_X)_j\Big((G_{t,\tau} \sigma Q_tX)_k\Big) X_k(G_{t,\tau}X)_j\\
= & \ \sum_{\substack{1 \leq j \leq n \\1 \leq k \leq 2n}}\partial_{\xi_j}\Big((G_{t,\tau} \sigma Q_tX)_k\Big) X_k(G_{t,\tau}X)_{j}
-\sum_{\substack{1 \leq j \leq n \\1 \leq k \leq 2n}}\partial_{x_j}\Big((G_{t,\tau} \sigma Q_tX)_k\Big) X_k(G_{t,\tau}X)_{j+n}\\
= & \ \sum_{\substack{1 \leq j \leq n \\1 \leq k \leq 2n}}(G_{t,\tau} \sigma Q_{t})_{k,j+n} X_k(G_{t,\tau}X)_{j}
-\sum_{\substack{1 \leq j \leq n \\1 \leq k \leq 2n}}(G_{t,\tau} \sigma Q_t)_{k,j} X_k(G_{t,\tau}X)_{j+n}\\
=& -\langle G_{t,\tau} \sigma Q_t \sigma G_{t,\tau}X,X\rangle.
\end{align*}
On the other hand, we observe that 
\begin{multline*}
\sum_{1 \leq j \leq 2n}(\sigma \nabla_X)_j\Big((G_{t,\tau} \sigma Q_tX)_j\Big)=\sum_{1 \leq j \leq n}\partial_{\xi_j}\Big((G_{t,\tau} \sigma Q_tX)_j\Big)-\sum_{1 \leq j \leq n}\partial_{x_j}\Big((G_{t,\tau} \sigma Q_tX)_{j+n}\Big)\\
=\sum_{1 \leq j \leq n}(G_{t,\tau} \sigma Q_t)_{j,j+n}-\sum_{1 \leq j \leq n}(G_{t,\tau} \sigma Q_t)_{j+n,j}=-\textrm{Tr}(\sigma G_{t,\tau} \sigma Q_t).
\end{multline*}
By using that the matrices $Q_t$ and $G_{t,\tau}$ are symmetric and $\sigma$ is skew-symmetric, the equations (\ref{jk8}), (\ref{jk8.k}), (\ref{jk9}) and (\ref{jk9.k}) reduce to 
\begin{multline}\label{jk10}
\partial_tG_{t,\tau}=Q_t-i\big(G_{t,\tau}\sigma Q_t+(G_{t,\tau}\sigma Q_t)^T\big)\\
+\frac{1}{2}\big(G_{t,\tau} \sigma Q_t \sigma G_{t,\tau}+(G_{t,\tau} \sigma Q_t \sigma G_{t,\tau})^T\big)
=Q_t-i(G_{t,\tau}\sigma Q_t-Q_t\sigma G_{t,\tau})+G_{t,\tau} \sigma Q_t \sigma G_{t,\tau},
\end{multline}
\begin{multline}\label{jk10.k}
\partial_{\tau}G_{t,\tau}=-Q_{\tau}-i\big(G_{t,\tau}\sigma Q_{\tau}+(G_{t,\tau}\sigma Q_{\tau})^T\big)\\
-\frac{1}{2}\big(G_{t,\tau} \sigma Q_{\tau} \sigma G_{t,\tau}+(G_{t,\tau} \sigma Q_{\tau} \sigma G_{t,\tau})^T\big)
=-Q_{\tau}-i(G_{t,\tau}\sigma Q_{\tau}-Q_{\tau}\sigma G_{t,\tau})-G_{t,\tau} \sigma Q_{\tau} \sigma G_{t,\tau},
\end{multline}
\begin{equation}\label{jk11}
\partial_th(t,\tau)=\frac{1}{2}\textrm{Tr}(\sigma G_{t,\tau} \sigma Q_t),
\end{equation}
\begin{equation}\label{jk11.k}
\partial_{\tau}h(t,\tau)=-\frac{1}{2}\textrm{Tr}(\sigma G_{t,\tau} \sigma Q_{\tau}),
\end{equation}
where $A^T$ denotes the transpose matrix of $A$.
By denoting $\tilde{S}_{t,\tau}=\sigma G_{t,\tau}$ the Hamilton map of the quadratic form $X \mapsto \langle G_{t,\tau} X,X\rangle$ and $F_t=\sigma Q_t$ the Hamilton map of the quadratic form $q_t(X)=\langle Q_t X,X\rangle$, we deduce from (\ref{jk10}), (\ref{jk10.k}), (\ref{jk11}) and (\ref{jk11.k}) that  
\begin{equation}\label{jk12}
\partial_t\tilde{S}_{t,\tau}=F_t-i(\tilde{S}_{t,\tau} F_t-F_t \tilde{S}_{t,\tau})+\tilde{S}_{t,\tau} F_t \tilde{S}_{t,\tau},
\end{equation}
\begin{equation}\label{jk12.k}
\partial_{\tau}\tilde{S}_{t,\tau}=-F_{\tau}-i(\tilde{S}_{t,\tau} F_{\tau}-F_{\tau} \tilde{S}_{t,\tau})-\tilde{S}_{t,\tau} F_{\tau} \tilde{S}_{t,\tau},
\end{equation}
\begin{equation}\label{jk13}
\partial_th(t,\tau)=\frac{1}{2}\textrm{Tr}(\tilde{S}_{t,\tau} F_t),
\end{equation}
\begin{equation}\label{jk13.k}
\partial_{\tau}h(t,\tau)=-\frac{1}{2}\textrm{Tr}(\tilde{S}_{t,\tau} F_{\tau}).
\end{equation}
We observe that the Hamilton map $\tilde{S}_{t,\tau}$ satisfies a matrix Ricatti differential equation. In order to solve this differential equation, we follow~\cite{riccati} (Chapter~2) and consider the first order linear differential equation
\begin{equation}\label{jk14}
Y'(t)=M(t)Y(t), \qquad Y(t)=\left(\begin{array}{c}Y_1(t) \\ Y_2(t)\end{array}\right) \in \cc^{4n \times 2n},
\end{equation}
with 
\begin{equation}\label{jk15}
M(t)=\left(\begin{array}{cc}
iF_t & -F_t\\
F_t  & i F_t
\end{array}\right) \in \cc^{4n \times 4n}.
\end{equation}
We observe that 
\begin{equation}\label{jk16}
\frac{d}{dt}\big(Y_1(t)-iY_2(t)\big)=0, \qquad \frac{d}{dt}\big(Y_1(t)+iY_2(t)\big)=2iF_t\big(Y_1(t)+iY_2(t)\big).
\end{equation}
With $R$ the resolvent
\begin{equation}\label{jk0}
\left\lbrace 
\begin{array}{ll}
\partial_tR(t,\tau)=2iF_tR(t,\tau), \qquad 0 \leq t \leq T,\\
R(\tau,\tau)=I_{2n},
\end{array}
\right.
\end{equation} 
with $0 \leq \tau \leq T$, we have
$$\forall 0 \leq t \leq T, \quad Y_1(t)-iY_2(t)=Y_1(\tau)-iY_2(\tau), \qquad Y_1(t)+iY_2(t)=R(t,\tau)(Y_1(\tau)+iY_2(\tau)).$$
It follows that 
$$\forall t \in [0,T], \quad Y_1(t)=\frac{1}{2}\big(R(t,\tau)+I_{2n}\big)Y_1(\tau)+\frac{i}{2}\big(R(t,\tau)-I_{2n}\big)Y_2(\tau),$$
$$\forall t \in [0,T], \quad Y_2(t)=\frac{1}{2i}\big(R(t,\tau)-I_{2n}\big)Y_1(\tau)+\frac{1}{2}\big(R(t,\tau)+I_{2n}\big)Y_2(\tau).$$
With the initial conditions $Y_1(\tau)=0$ and $Y_2(\tau)=I_{2n}$, this leads to consider the function
\begin{equation}\label{jk17}
S(t,\tau)=-Y_1(t)Y_2(t)^{-1}=-i\big(R(t,\tau)-I_{2n}\big)\big(R(t,\tau)+I_{2n}\big)^{-1},
\end{equation}
which is well-defined when $|t-\tau| \ll 1$ is sufficiently small, since $R(\tau,\tau)=I_{2n}$.
By differentiating the identity 
$$I_{2n}=\big(R(t,\tau)+I_{2n}\big)^{-1}\big(R(t,\tau)+I_{2n}\big),$$
we obtain that 
\begin{equation}\label{jk18}
\frac{d}{dt}\big(R(t,\tau)+I_{2n}\big)^{-1}=-2i\big(R(t,\tau)+I_{2n}\big)^{-1}F_tR(t,\tau)\big(R(t,\tau)+I_{2n}\big)^{-1},
\end{equation}
when $|t-\tau| \ll 1$.
It follows from (\ref{jk17}) and (\ref{jk18}) that 
\begin{align*}
\partial_tS(t,\tau)=& \ 2F_tR(t,\tau)\big(R(t,\tau)+I_{2n}\big)^{-1}\\
& \qquad -2\big(R(t,\tau)-I_{2n}\big)\big(R(t,\tau)+I_{2n}\big)^{-1}F_tR(t,\tau)\big(R(t,\tau)+I_{2n}\big)^{-1}\\
= & \ 4\big(R(t,\tau)+I_{2n}\big)^{-1}F_tR(t,\tau)\big(R(t,\tau)+I_{2n}\big)^{-1}\\
= & \ 4\big(R(t,\tau)+I_{2n}\big)^{-1}F_t-4\big(R(t,\tau)+I_{2n}\big)^{-1}F_t\big(R(t,\tau)+I_{2n}\big)^{-1},
\end{align*}
when $|t-\tau| \ll 1$.
On the other hand, we deduce from (\ref{jk17}) that 
\begin{align*}
& \  F_t-i\big(S(t,\tau)F_t-F_t S(t,\tau)\big)+S(t,\tau)F_t S(t,\tau)\\
= & \ F_t-\big(R(t,\tau)-I_{2n}\big)\big(R(t,\tau)+I_{2n}\big)^{-1}F_t
 +F_t\big(R(t,\tau)-I_{2n}\big)\big(R(t,\tau)+I_{2n}\big)^{-1}\\
&\ -\big(R(t,\tau)-I_{2n}\big)\big(R(t,\tau)+I_{2n}\big)^{-1}F_t\big(R(t,\tau)-I_{2n}\big)\big(R(t,\tau)+I_{2n}\big)^{-1},
\end{align*}
when $|t-\tau| \ll 1$.
A direct computation provides
\begin{align*}
& \  F_t-i\big(S(t,\tau)F_t-F_t S(t,\tau)\big)+S(t,\tau)F_t S(t,\tau)\\
= &\  2\big(R(t,\tau)+I_{2n}\big)^{-1}F_t
 +F_t-2F_t\big(R(t,\tau)+I_{2n}\big)^{-1} -F_t\big(R(t,\tau)-I_{2n}\big)\big(R(t,\tau)+I_{2n}\big)^{-1}\\
 & \ +2\big(R(t,\tau)+I_{2n}\big)^{-1}F_t\big(R(t,\tau)-I_{2n}\big)\big(R(t,\tau)+I_{2n}\big)^{-1},
\end{align*}
implying that 
\begin{align*}
& \  F_t-i\big(S(t,\tau)F_t-F_t S(t,\tau)\big)+S(t,\tau)F_t S(t,\tau)\\
= &\  2\big(R(t,\tau)+I_{2n}\big)^{-1}F_t+2\big(R(t,\tau)+I_{2n}\big)^{-1}F_t\big(R(t,\tau)-I_{2n}\big)\big(R(t,\tau)+I_{2n}\big)^{-1}\\
= &\  4\big(R(t,\tau)+I_{2n}\big)^{-1}F_t-4\big(R(t,\tau)+I_{2n}\big)^{-1}F_t\big(R(t,\tau)+I_{2n}\big)^{-1},
\end{align*}
when $|t-\tau| \ll 1$.
We therefore notice that the function $t \mapsto S(t,\tau)$ defined in (\ref{jk17}) satisfies the differential equation (\ref{jk12}).
On the other hand, we recall for instance from~\cite{coron} (Proposition~1.5) that the resolvent satisfies 
\begin{equation}\label{kl1}
\partial_{\tau}R(t,\tau)=-2iR(t,\tau)F_{\tau}.
\end{equation}
By differentiating the identity 
$$I_{2n}=\big(R(t,\tau)+I_{2n}\big)^{-1}\big(R(t,\tau)+I_{2n}\big),$$
we obtain that 
\begin{equation}\label{jk18.k}
\frac{d}{d\tau}\big(R(t,\tau)+I_{2n}\big)^{-1}=2i\big(R(t,\tau)+I_{2n}\big)^{-1}R(t,\tau)F_{\tau}\big(R(t,\tau)+I_{2n}\big)^{-1},
\end{equation}
when $|t-\tau| \ll 1$.
It follows from (\ref{jk17}) and (\ref{jk18.k}) that 
\begin{align*}
\partial_{\tau}S(t,\tau)=& \ -2R(t,\tau)F_{\tau}\big(R(t,\tau)+I_{2n}\big)^{-1}\\
& \qquad +2\big(R(t,\tau)-I_{2n}\big)\big(R(t,\tau)+I_{2n}\big)^{-1}R(t,\tau)F_{\tau}\big(R(t,\tau)+I_{2n}\big)^{-1}\\
= & \ -4\big(R(t,\tau)+I_{2n}\big)^{-1}R(t,\tau)F_{\tau}\big(R(t,\tau)+I_{2n}\big)^{-1}\\
= & \ -4F_{\tau}\big(R(t,\tau)+I_{2n}\big)^{-1}+4\big(R(t,\tau)+I_{2n}\big)^{-1}F_{\tau}\big(R(t,\tau)+I_{2n}\big)^{-1},
\end{align*}
when $|t-\tau| \ll 1$.
On the other hand, we deduce from (\ref{jk17}) that 
\begin{align*}
& \  -F_{\tau}-i\big(S(t,\tau)F_{\tau}-F_{\tau} S(t,\tau)\big)-S(t,\tau)F_{\tau} S(t,\tau)\\
= & \ -F_{\tau}-\big(R(t,\tau)-I_{2n}\big)\big(R(t,\tau)+I_{2n}\big)^{-1}F_{\tau}
 +F_{\tau}\big(R(t,\tau)-I_{2n}\big)\big(R(t,\tau)+I_{2n}\big)^{-1}\\
&\ +\big(R(t,\tau)-I_{2n}\big)\big(R(t,\tau)+I_{2n}\big)^{-1}F_{\tau}\big(R(t,\tau)-I_{2n}\big)\big(R(t,\tau)+I_{2n}\big)^{-1},
\end{align*}
when $|t-\tau| \ll 1$.
A direct computation provides
\begin{align*}
& \  -F_{\tau}-i\big(S(t,\tau)F_{\tau}-F_{\tau} S(t,\tau)\big)-S(t,\tau)F_{\tau} S(t,\tau)\\
= &\  2\big(R(t,\tau)+I_{2n}\big)^{-1}F_{\tau}
 -2F_{\tau}\big(R(t,\tau)+I_{2n}\big)^{-1} +\big(R(t,\tau)-I_{2n}\big)\big(R(t,\tau)+I_{2n}\big)^{-1}F_{\tau}\\
 & \ -F_{\tau} -2\big(R(t,\tau)-I_{2n}\big)\big(R(t,\tau)+I_{2n}\big)^{-1}F_{\tau}\big(R(t,\tau)+I_{2n}\big)^{-1},
\end{align*}
implying that 
\begin{align*}
& \  -F_{\tau}-i\big(S(t,\tau)F_{\tau}-F_{\tau} S(t,\tau)\big)-S(t,\tau)F_{\tau} S(t,\tau)\\
= &\  -2F_{\tau}\big(R(t,\tau)+I_{2n}\big)^{-1}-2\big(R(t,\tau)-I_{2n}\big)\big(R(t,\tau)+I_{2n}\big)^{-1}F_{\tau}\big(R(t,\tau)+I_{2n}\big)^{-1}\\
= &\  -4F_{\tau}\big(R(t,\tau)+I_{2n}\big)^{-1}+4\big(R(t,\tau)+I_{2n}\big)^{-1}F_{\tau}\big(R(t,\tau)+I_{2n}\big)^{-1},
\end{align*}
when $|t-\tau| \ll 1$.
We therefore notice that the function $\tau \mapsto S(t,\tau)$ defined in (\ref{jk17}) satisfies the differential equation (\ref{jk12.k}).
Let $\textrm{Log }z$ be the principal determination of the complex logarithm on $\cc \setminus \rr_-$. We consider the function
\begin{equation}\label{jk19}
h(t,\tau)=-\frac{1}{2}\textrm{Log}\big(2^{-2n}\textrm{det}\big(R(t,\tau)+I_{2n}\big)\big),
\end{equation} 
which is well-defined when $|t-\tau| \ll 1$, since $R(\tau,\tau)=I_{2n}$. With $\textrm{Com}(A)$ denoting the adjugate matrix of $A$, that is, the transpose of the cofactor matrix of $A$, we indeed notice  from (\ref{jk17}) that it satisfies
\begin{align*}
& \ \partial_th(t,\tau)= -\frac{1}{2}\big(\textrm{det}\big(R(t,\tau)+I_{2n}\big)\big)^{-1}\textrm{Tr}\big(\big[\textrm{Com}\big(R(t,\tau)+I_{2n}\big)\big]^T(2i)F_tR(t,\tau)\big)\\
& \ =-i\textrm{Tr}\big(\big(R(t,\tau)+I_{2n}\big)^{-1}F_tR(t,\tau)\big)=\frac{1}{2}\textrm{Tr}\big(F_tS(t,\tau)\big)-\frac{i}{2}\textrm{Tr}\big(F_t)=\frac{1}{2}\textrm{Tr}\big(S(t,\tau)F_t\big),
\end{align*} 
when $|t-\tau| \ll 1$, since from (\ref{vb1}), we have $\textrm{Tr}(F_t)=0$. It proves the formula (\ref{jk13}). On the other hand, we deduce from (\ref{kl1}) that  
\begin{align*}
 \partial_{\tau}h(t,\tau)= & \ \frac{1}{2}\big(\textrm{det}\big(R(t,\tau)+I_{2n}\big)\big)^{-1}\textrm{Tr}\big(\big[\textrm{Com}\big(R(t,\tau)+I_{2n}\big)\big]^T(2i)R(t,\tau)F_{\tau}\big)\\
 = & \ i\textrm{Tr}\big(\big(R(t,\tau)+I_{2n}\big)^{-1}R(t,\tau)F_{\tau}\big)\\
 = & \ \frac{i}{2}\textrm{Tr}(F_{\tau})+\frac{i}{2}\textrm{Tr}\big(\big(R(t,\tau)+I_{2n}\big)^{-1}\big(R(t,\tau)-I_{2n}\big)F_{\tau}\big),
\end{align*} 
when $|t-\tau| \ll 1$.
We notice 
\begin{align}\label{jk21}
& \ S(t,\tau)=-i\big(R(t,\tau)-I_{2n}\big)\big(R(t,\tau)+I_{2n}\big)^{-1}\\ \notag
=& \ i\frac{I_{2n}-R(t,\tau)}{2}\Big(I_{2n}-\frac{I_{2n}-R(t,\tau)}{2}\Big)^{-1}
=i\sum_{k=0}^{+\infty}\frac{1}{2^{k+1}}\big(I_{2n}-R(t,\tau)\big)^{k+1}\\ \notag
=& \  i\Big(I_{2n}-\frac{I_{2n}-R(t,\tau)}{2}\Big)^{-1}\frac{I_{2n}-R(t,\tau)}{2}=-i\big(R(t,\tau)+I_{2n}\big)^{-1}\big(R(t,\tau)-I_{2n}\big),
\end{align}
when $|t-\tau| \ll 1$, since $R(\tau,\tau)=I_{2n}$. It follows from (\ref{jk21}) that
\begin{align*}
& \ \partial_{\tau}h(t,\tau)= -\frac{1}{2}\textrm{Tr}\big(S(t,\tau)F_{\tau}\big),
\end{align*} 
since $\textrm{Tr}(F_{\tau})=0$. It proves the formula (\ref{jk13.k}).
We need the following instrumental lemma:

\bigskip

\begin{lemma}\label{lem1}
Let $R(t,\tau)$ be the resolvent 
$$\left\lbrace 
\begin{array}{ll}
\frac{d}{dt}R(t,\tau)=2iF_tR(t,\tau), \qquad 0 \leq t \leq T,\\
R(\tau,\tau)=I_{2n},
\end{array}
\right.
$$
with $0 \leq \tau \leq T$.
Then, the mapping $R(t,\tau) : \cc^{2n} \rightarrow \cc^{2n}$ is a non-negative complex symplectic linear transformation satisfying
$$\forall t,\tau \in [0,T],  \quad R(t,\tau)^{-1}=R(\tau,t),$$
$$\forall t,\tau \in [0,T], \forall X,Y \in \cc^{2n}, \quad \sigma(R(t,\tau)X,R(t,\tau)Y)=\sigma(X,Y),$$
$$\forall 0 \leq \tau \leq t \leq T, \forall X \in \cc^{2n}, \quad i\big(\sigma(\overline{R(t,\tau)X},R(t,\tau)X)-\sigma(\overline{X},X)\big) \geq 0.$$
\end{lemma}

\bigskip

\begin{proof}
Standard results about resolvents show that the mapping $R(t,\tau) : \cc^{2n} \rightarrow \cc^{2n}$ defines an isomorphism whose inverse is $R(t,\tau)^{-1}=R(\tau,t)$. On the other hand, we notice from (\ref{a1}) and (\ref{jk0}) that for all $0 \leq t, \tau \leq T$,
\begin{multline*}
\frac{d}{dt}\big(\sigma(R(t,\tau)X,R(t,\tau)Y)\big)=\sigma(2iF_tR(t,\tau)X,R(t,\tau)Y)+\sigma(R(t,\tau)X,2iF_tR(t,\tau)Y)\\
=2i\sigma(F_tR(t,\tau)X,R(t,\tau)Y)-2i\sigma(F_tR(t,\tau)X,R(t,\tau)Y)=0.
\end{multline*}
By using that $\sigma(R(\tau,\tau)X,R(\tau,\tau)Y)=\sigma(X,Y),$ since $R(\tau,\tau)=I_{2n}$, we obtain that 
$$\forall 0 \leq t, \tau \leq T, \forall X,Y \in \cc^{2n}, \quad \sigma(R(t,\tau)X,R(t,\tau)Y)=\sigma(X,Y).$$
Setting 
$$f_{\tau}(t)=i\big(\sigma(\overline{R(t,\tau)X},R(t,\tau)X)-\sigma(\overline{X},X)\big), \quad 0 \leq t \leq T,\ X \in \cc^{2n},$$
with $0 \leq \tau \leq T$, we observe that $f_{\tau}(\tau)=0$, since $R(\tau,\tau)=I_{2n}$. On the other hand, it follows from (\ref{a1}) and (\ref{jk0}) that for all $0 \leq t \leq T$,
\begin{align*}
& \ f_{\tau}'(t)=i\sigma(\overline{2iF_tR(t,\tau)X},R(t,\tau)X)+i\sigma(\overline{R(t,\tau)X},2iF_tR(t,\tau)X)\\
= &\ -2\sigma(\overline{R(t,\tau)X},(F_t+\overline{F_t})R(t,\tau)X)=-4\sigma(\overline{R(t,\tau)X},\textrm{Re }F_tR(t,\tau)X)\\
=& \ -4(\textrm{Re }q_t)(\overline{R(t,\tau)X},R(t,\tau)X) =-4\langle \textrm{Re }Q_t \overline{R(t,\tau)X},R(t,\tau)X\rangle\\
= &\ -4\langle \textrm{Re }Q_t \textrm{Re}(R(t,\tau)X),\textrm{Re}(R(t,\tau)X)\rangle-4\langle \textrm{Re }Q_t \textrm{Im}(R(t,\tau)X),\textrm{Im}(R(t,\tau)X)\rangle \geq 0,
\end{align*}
since $\textrm{Re }Q_t \leq 0$.
We deduce that 
$$\forall 0 \leq \tau \leq t \leq T, \quad f_{\tau}(t)=i\big(\sigma(\overline{R(t,\tau)X},R(t,\tau)X)-\sigma(\overline{X},X)\big) \geq 0.$$
This ends the proof of Lemma~\ref{lem1}.
\end{proof}

\noindent
The following lemma shows that the matrix 
$$S(t,\tau)=-i\big(R(t,\tau)-I_{2n}\big)\big(R(t,\tau)+I_{2n}\big)^{-1},$$ 
defined in (\ref{jk17}) is a Hamilton map:

\bigskip

\begin{lemma}\label{lem2}
The matrix 
$$S(t,\tau)=-i\big(R(t,\tau)-I_{2n}\big)\big(R(t,\tau)+I_{2n}\big)^{-1},$$ 
defined for all $0 \leq \tau \leq t \leq T$ and $0 \leq t-\tau \leq \delta$, with $0<\delta \ll 1$, is the Hamilton map associated to the quadratic form 
$$X \in \rr^{2n} \mapsto \langle G_{t,\tau}X,X\rangle=\sigma(X,S(t,\tau)X\big) \in \cc,$$
whose real part is non-positive
$$\forall 0 \leq \tau \leq t \leq T, \forall X \in \rr^{2n}, \quad  \emph{\textrm{Re}}(\langle G_{t,\tau}X,X\rangle) \leq 0.$$
\end{lemma}

\bigskip

\begin{proof}
It follows from Lemma~\ref{lem1} that 
\begin{equation}\label{jk22}
\forall 0 \leq t,\tau \leq T, \forall X,Y \in \cc^{2n}, \quad \sigma(R(t,\tau)X,Y)=\sigma(X,R(\tau,t)Y).
\end{equation}
We deduce from (\ref{jk21}) and (\ref{jk22}) that 
\begin{equation}\label{jk23}
\forall X,Y \in \cc^{2n}, \quad \sigma(S(t,\tau)X,Y)=\sigma(X,S(\tau,t)Y),
\end{equation}
when $|t-\tau| \ll 1$, since
$$S(t,\tau)=i\sum_{k=0}^{+\infty}\frac{1}{2^{k+1}}\big(I_{2n}-R(t,\tau)\big)^{k+1}.$$
We want to prove that the matrix $S(t,\tau)$ is the Hamilton map associated to the quadratic form 
$$X \mapsto \sigma\big(X,S(t,\tau)X\big).$$ 
According to (\ref{10}), (\ref{a1}) and (\ref{jk23}), it is sufficient to establish that $S(t,\tau)=-S(\tau,t)$, when $|t-\tau| \ll 1$.
By using (\ref{jk21}), this is equivalent to the following identity
$$-\big(R(t,\tau)-I_{2n}\big)\big(R(t,\tau)+I_{2n}\big)^{-1}=\big(R(\tau,t)+I_{2n}\big)^{-1}\big(R(\tau,t)-I_{2n}\big),$$
that is
$$-\big(R(\tau,t)+I_{2n}\big)\big(R(t,\tau)-I_{2n}\big)=\big(R(\tau,t)-I_{2n}\big)\big(R(t,\tau)+I_{2n}\big),$$
which holds true since
$$-\big(R(\tau,t)+I_{2n}\big)\big(R(t,\tau)-I_{2n}\big)=R(\tau,t)-R(t,\tau)=\big(R(\tau,t)-I_{2n}\big)\big(R(t,\tau)+I_{2n}\big),$$
since $R(t_1,t_2)R(t_2,t_3)=R(t_1,t_3)$ when $0 \leq t_1,t_2,t_3 \leq T$. On the other hand, we deduce from Lemma~\ref{lem1} that for all $X \in \cc^{2n}$, $0 \leq \tau \leq t \leq T$,
\begin{align*}
& \ \textrm{Re}\big(i\sigma\big(\overline{\big(R(t,\tau)+I_{2n}\big)X},\big(R(t,\tau)-I_{2n}\big)X\big)\big)\\
= & \ \textrm{Re}\big(i\big[\sigma\big(\overline{R(t,\tau)X},R(t,\tau)X\big)-\sigma\big(\overline{X},X\big)\big]\big)+\textrm{Re}\big(i\big[\sigma\big(\overline{X},R(t,\tau)X\big)-\sigma\big(\overline{R(t,\tau)X},X\big)\big]\big)\\
= & \ i\big[\sigma\big(\overline{R(t,\tau)X},R(t,\tau)X\big)-\sigma\big(\overline{X},X\big)\big]+\textrm{Re}\big(i\big[\sigma\big(\overline{X},R(t,\tau)X\big)+\overline{\sigma\big(\overline{X},R(t,\tau)X\big)}\big]\big)\\
= & \ i\big[\sigma\big(\overline{R(t,\tau)X},R(t,\tau)X\big)-\sigma\big(\overline{X},X\big)\big] \geq 0.
\end{align*}
We deduce from the above estimate that  
\begin{equation}\label{jk24}
\forall X \in \cc^{2n} , \quad \textrm{Re}\big(i\sigma\big(\overline{X},\big(R(t,\tau)-I_{2n}\big)\big(R(t,\tau)+I_{2n}\big)^{-1}X\big)\big) \geq 0,
\end{equation}
when $0 \leq \tau \leq t \leq T$ and $|t-\tau| \ll 1$. We obtain in particular from (\ref{jk24}) that 
$$\forall X \in \rr^{2n} , \quad \textrm{Re}\big(\sigma\big(X,S(t,\tau)X\big)\big) =\textrm{Re}\big(-i\sigma\big(X,\big(R(t,\tau)-I_{2n}\big)\big(R(t,\tau)+I_{2n}\big)^{-1}X\big)\big) \leq 0,$$
when $0 \leq \tau \leq t \leq T$ and $|t-\tau| \ll 1$. This ends the proof of Lemma~\ref{lem2}. 

\end{proof}

We consider the Weyl symbol 
\begin{equation}\label{weyl}
p_{t,\tau}(X)=\frac{2^n}{\sqrt{\textrm{det}\big(R(t,\tau)+I_{2n}\big)}}\exp\big(-i\sigma(X,\big(R(t,\tau)-I_{2n}\big)\big(R(t,\tau)+I_{2n}\big)^{-1}X\big)\big),
\end{equation}
with $X=(x,\xi) \in \rr^{2n}$, for all $0 \leq t,\tau \leq T$, $|t-\tau| \leq \delta$, where the positive constant $\delta >0$ is chosen sufficiently small for the determinant $\textrm{det}\big(R(t,\tau)+I_{2n}\big) \neq 0$ to be non-zero and its square root well-defined when using the principal determination of the complex logarithm. This is possible as $R(t,t)=I_{2n}$ when $0 \leq t \leq T$.
We notice from (\ref{jk2}), (\ref{jk3}), (\ref{jk19}) and Lemma~\ref{lem2} that it is equal to the symbol 
$$p_{t,\tau}(X)=e^{g_{t,\tau}(X)}=\exp\big(\langle G_{t,\tau}X,X\rangle+h(t,\tau)\big),$$ 
and therefore satisfies the equations
\begin{equation}\label{carre}
\frac{d}{dt}p_{t,\tau}=q_t \#^w p_{t,\tau}, \qquad \frac{d}{d\tau}p_{t,\tau}=-p_{t,\tau} \#^w q_{\tau} , \qquad p_{\tau,\tau}=1,
\end{equation}
when $0 \leq t,\tau \leq T$, $|t-\tau| \leq \delta$. On the other hand, notice that Lemma~\ref{lem2} implies that 
\begin{multline}\label{hj1}
\forall 0 \leq \tau \leq t \leq T, \ 0 \leq t-\tau \leq \delta, \forall X \in \rr^{2n}, \\ 
|\exp(\langle G_{t,\tau}X,X\rangle)|=\big|\exp\big(-i\sigma(X,\big(R(t,\tau)-I_{2n}\big)\big(R(t,\tau)+I_{2n}\big)^{-1}X\big)\big)\big| \leq 1.
\end{multline}
The symbol $p_{t,\tau}$ is therefore a $L^{\infty}(\rr_X^{2n})$-function when $0 \leq \tau \leq t \leq T$, $0 \leq t-\tau \leq \delta$.

We consider the pseudodifferential operator $p_{t,\tau}^w(x,D_x)$ defined by the Weyl quantization of the symbol $p_{t,\tau}$. 
We aim at proving that this pseudodifferential operator is equal to the Fourier integral operator 
$$\mathscr{K}_{R(\tau,t)} : \mathscr{S}(\rr^n) \rightarrow \mathscr{S}'(\rr^n),$$ associated to the non-negative complex symplectic linear transformation $R(t,\tau)$. Setting
$$\tilde{S}(t,\tau)=-\big(R(t,\tau)-I_{2n}\big)\big(R(t,\tau)+I_{2n}\big)^{-1},$$
the following identities
\begin{equation}\label{ll1}
I_{2n}+\tilde{S}(t,\tau)=2\big(R(t,\tau)+I_{2n}\big)^{-1} \quad \textrm{and} \quad  I_{2n}-\tilde{S}(t,\tau)=2R(t,\tau)\big(R(t,\tau)+I_{2n}\big)^{-1},
\end{equation}
imply that 
\begin{equation}\label{jk30}
\big(I_{2n}-\tilde{S}(t,\tau)\big)\big(I_{2n}+\tilde{S}(t,\tau)\big)^{-1}=R(t,\tau),
\end{equation}
when $0 \leq \tau \leq t \leq T$, $0 \leq t-\tau \leq \delta$.
We observe that $\tilde{S}(\tau,\tau)=0$, since $R(\tau,\tau)=I_{2n}$ for $0 \leq \tau \leq T$. By possibly decreasing the value of the positive constant $\delta>0$, it follows that $\pm 1$ are not eigenvalues of the matrix $\tilde{S}(t,\tau)$ for all $0 \leq \tau \leq t \leq T$, $0 \leq t-\tau \leq \delta$. 
We can therefore deduce from the link between pseudodifferential operators and Fourier integral operators established by H\"ormander in~\cite{mehler} (Proposition~5.11), Lemma~\ref{lem2}, (\ref{weyl}),  (\ref{ll1}) and (\ref{jk30}) that for all $0 \leq \tau \leq t \leq T$, $0 \leq t-\tau \leq \delta$, 
\begin{multline}\label{jk31}
\mathscr{K}_{R(\tau,t)}=\sqrt{\frac{2^{2n}\textrm{det}\big(R(t,\tau)\big)}{\textrm{det}\big(R(t,\tau)+I_{2n}\big)}}\big(e^{-i\sigma(X,(R(t,\tau)-I_{2n})(R(t,\tau)+I_{2n})^{-1}X)}\big)^w\\
=\frac{2^{n}}{\sqrt{\textrm{det}\big(R(t,\tau)+I_{2n}\big)}}\big(e^{-i\sigma(X,(R(t,\tau)-I_{2n})(R(t,\tau)+I_{2n})^{-1}X)}\big)^w=p_{t,\tau}^w(x,D_x),
\end{multline}
since $\textrm{det}\big(R(t,\tau)\big)=1$, because $R(t,\tau) : \cc^{2n} \rightarrow \cc^{2n}$ is a non-negative complex symplectic linear transformation and therefore belongs to the special linear group $\textrm{SL}_{2n}(\cc)$. Indeed, the real symplectic linear group is included in the real special linear group $\textrm{SL}_{2n}(\rr)$, see e.g.~\cite{Le} (Proposition~4.4.4). On the other hand, we know from~\cite{mehler} (Proposition~5.10) that any non-negative complex symplectic linear transformation $\mathcal{T} : \cc^{2n} \rightarrow \cc^{2n}$ can be factored as $\mathcal{T}=\mathcal{T}_1\mathcal{T}_2\mathcal{T}_3$, where $\mathcal{T}_1$ and $\mathcal{T}_3$ are real symplectic linear transformations and $\mathcal{T}_2(x,\xi)=(x',\xi')$ where for all $1 \leq j \leq n$, either
$$(x_j',\xi_j')=(x_j \cosh \tau_j-i\xi_j \sinh \tau_j, ix_j \sinh \tau_j+\xi_j \cosh \tau_j),$$
with $\tau_j \geq 0$, or 
$$(x_j',\xi_j')=(x_j,ix_j+\xi_j).$$

We consider 
$$\chi_{\eps}(x,\xi)=\chi(\eps x,\eps \xi),$$ 
where $\chi \in C_0^{\infty}(\rr^{2n},\rr )$ is equal to $1$ in a neighborhood of $0$.
By Calder\'on-Vaillancourt Theorem, the pseudodifferential operator $\chi_{\eps}^w(x,D_x)$ defines a bounded selfadjoint operator on $L^2(\rr^n)$, whose operator norm is uniformly bounded with respect to the parameter $0<\eps \leq 1$,
\begin{equation}\label{as2}
\exists C>0, \forall 0<\eps \leq 1, \quad \|\chi_{\eps}^w(x,D_x)\|_{\mathcal{L}(L^2)} \leq C.
\end{equation}
Furthermore, it is also a continuous mapping from $L^2(\rr^n)$ to $\mathscr{S}(\rr^n)$ since $\chi_{\eps} \in \mathscr{S}(\rr^{2n})$.
We observe that the symbol $(\chi_{\eps})_{0<\eps \leq 1}$ is bounded in the Fr\'echet space $C_b^{\infty}(\rr^{2n})$ and that $(\chi_{\eps})_{0<\eps \leq 1}$ converges in $C^{\infty}(\rr^{2n})$ to the constant function $1$, when $\eps$ tends to $0$. It follows from~\cite{Le} (Lemma~1.1.3) that the sequence $(\chi_{\eps}^w(x,D_x)u)_{0<\eps \leq 1}$ converges to $u$ in $\mathscr{S}(\rr^n)$, if $u \in \mathscr{S}(\rr^n)$. On the other hand, it follows from (\ref{as2}) that for all $u \in L^2(\rr^n)$ and $v \in \mathscr{S}(\rr^n)$, 
\begin{multline}\label{conv}
\limsup_{\eps \to 0}\|u-\chi_{\eps}^w(x,D_x)u\|_{L^2(\rr^n)} \\
\leq \limsup_{\eps \to 0}\|v-\chi_{\eps}^w(x,D_x)v\|_{L^2(\rr^n)}+(C+1)\|u-v\|_{L^2(\rr^n)}  \leq (C+1)\|u-v\|_{L^2(\rr^n)}.
\end{multline}
By density of the Schwartz space in $L^2(\rr^n)$, we obtain that when $u \in L^2(\rr^n)$, the sequence $(\chi_{\eps}^w(x,D_x)u)_{0<\eps \leq 1}$ converges to $u$ in $L^2(\rr^n)$ when $\eps$ tends to $0$,
\begin{equation}\label{tria}
\forall u \in L^2(\rr^n), \quad \lim_{\eps \to 0}\|\chi_{\eps}^w(x,D_x)u-u\|_{L^2(\rr^n)}=0.
\end{equation}

Let $u, v \in \mathscr{S}(\rr^n)$. We deduce from Proposition~\ref{FIOdef} and (\ref{jk31}) that the function $p_{t,\tau}^w(x,D_x)u$ belongs to the Schwartz space for all $0 \leq \tau \leq t \leq T$, $0 \leq t-\tau \leq \delta$.
The theorem of regularity of integrals with parameters allows to obtain that for all $0 \leq \tau \leq t \leq T$, $0 \leq t-\tau \leq \delta$,
\begin{multline}\label{as4}
\frac{d}{d\tau}(p_{t,\tau}^w(x,D_x)u,v)_{L^2(\rr^n)} =\frac{d}{d\tau}\langle p_{t,\tau}^w(x,D_x)u,\overline{v} \rangle_{\mathscr{S}'(\rr^n),\mathscr{S}(\rr^n)}\\
=\frac{d}{d\tau} \int_{\rr^{2n}}p_{t,\tau}(x,\xi)\mathcal{H}(u,v)(x,\xi)dxd\xi= \int_{\rr^{2n}}\frac{d}{d\tau}p_{t,\tau}(x,\xi)\mathcal{H}(u,v)(x,\xi)dxd\xi,
\end{multline}
where $\mathcal{H}(u,v)$ denotes the Wigner function which defines a continuous mapping
\begin{equation}\label{wigner}
(u,v) \in \mathscr{S}(\rr^n) \times \mathscr{S}(\rr^n) \mapsto  \mathcal{H}(u,v) \in \mathscr{S}(\rr^{2n}),
\end{equation}
between the Schwartz spaces, see e.g.~\cite{Le} (Chapter~2). The differentiation under the integral sign in (\ref{as4}) is then justified as we notice from  (\ref{weyl}) and (\ref{hj1}) that
\begin{multline}\label{ku1}
\exists C_0>0, \forall 0 \leq \tau \leq t \leq T,\ 0 \leq t-\tau \leq \delta, \forall (x,\xi) \in \rr^{2n}, \\  \Big|\frac{d}{d\tau}p_{t,\tau}(x,\xi)\Big| \leq C_0(1+|x|^2+|\xi|^2).
\end{multline}
For $u, v \in \mathscr{S}(\rr^n)$, we define the function
\begin{equation}\label{as1}
f_{\eps}(\tau)=\big(p_{t_0,\tau}^w(x,D_x)\chi_{\eps}^w(x,D_x)\mathscr{U}(\tau,\tau_0)u,v\big)_{L^2(\rr^n)},
\end{equation}
when $\tau_0 \leq \tau \leq t_0$, with  $0 \leq \tau_0 < t_0 \leq T$, $0<t_0-\tau_0 \leq \delta$, where $(\mathscr{U}(t,\tau))_{0 \leq \tau \leq t \leq T}$ stands for the contraction evolution system given by Theorem~\ref{th1}.
This function is well-defined since $\mathscr{U}(\tau,\tau_0)u \in L^2(\rr^n)$ implies that 
\begin{equation}\label{hj2}
\chi_{\eps}^w(x,D_x)\mathscr{U}(\tau,\tau_0)u \in \mathscr{S}(\rr^n),
\end{equation} 
and, as Proposition~\ref{FIOdef} and (\ref{jk31}) provide that 
$$\forall \tau_0 \leq \tau \leq t_0, \quad p_{t_0,\tau}^w(x,D_x)\chi_{\eps}^w(x,D_x)\mathscr{U}(\tau,\tau_0)u \in \mathscr{S}(\rr^n).$$ 
We observe from (\ref{weyl}), (\ref{hj1}), (\ref{as1}) and (\ref{hj2}) that the mapping 
\begin{multline}\label{as6}
f_{\eps}(\tau)=\big(p_{t_0,\tau}^w(x,D_x)\chi_{\eps}^w(x,D_x)\mathscr{U}(\tau,\tau_0)u,v\big)_{L^2(\rr^n)}\\
=\int_{\rr^{2n}}p_{t_0,\tau}(x,\xi)\mathcal{H}(\chi_{\eps}^w(x,D_x)\mathscr{U}(\tau,\tau_0)u,v)(x,\xi)dxd\xi,
\end{multline}
is continuous on $[\tau_0,t_0]$. Indeed, we notice from (\ref{wigner}) that $\mathcal{H}(\chi_{\eps}^w(x,D_x)\mathscr{U}(\tau,\tau_0)u,v) \in \mathscr{S}(\rr^{2n})$ since $\chi_{\eps}^w(x,D_x)\mathscr{U}(\tau,\tau_0)u \in \mathscr{S}(\rr^n)$ and $v \in \mathscr{S}(\rr^n)$.
Furthermore, we deduce anew from (\ref{wigner}) that the continuity of the mapping $\tau \mapsto \mathscr{U}(\tau,\tau_0)u \in L^2(\rr^n)$ successively implies the continuity of the mappings $\tau \mapsto \chi_{\eps}^w(x,D_x)\mathscr{U}(\tau,\tau_0)u \in \mathscr{S}(\rr^n)$ and $\tau \mapsto \mathcal{H}(\chi_{\eps}^w(x,D_x)\mathscr{U}(\tau,\tau_0)u,v) \in \mathscr{S}(\rr^{2n})$. The domination condition then easily follows from the fact that any Schwartz seminorm of the Wigner function can be bounded as 
\begin{multline*}
\sup_{\substack{x,\xi \in \rr^n, \\ |\alpha_1|+|\alpha_2|+|\beta_1|+|\beta_2| \leq N_1}}|x^{\alpha_1}\xi^{\alpha_2}\partial_{x}^{\beta_1}\partial_{\xi}^{\beta_2}\mathcal{H}(\chi_{\eps}^w(x,D_x)\mathscr{U}(\tau,\tau_0)u,v)(x,\xi)| \\
\leq c \|\mathscr{U}(\tau,\tau_0)u\|_{L^2(\rr^n)}\Big(\sup_{\substack{x \in \rr^n, \\ |\alpha|+|\beta| \leq N_2}}|x^{\alpha}\partial_x^{\beta}v(x)|\Big)
\leq c \|u\|_{L^2(\rr^n)}\Big(\sup_{\substack{x \in \rr^n, \\ |\alpha|+|\beta| \leq N_2}}|x^{\alpha}\partial_x^{\beta}v(x)|\Big),
\end{multline*}
since $\|\mathscr{U}(\tau,\tau_0)\|_{\mathcal{L}(L^2)} \leq 1$.
On the other hand, we have for all $\tau_0<\tau<t_0$ and $0 \neq |h| \leq  \inf(t_0-\tau,\tau-\tau_0)$,
\begin{align}\label{as9}
& \ \frac{f_{\eps}(\tau+h)-f_{\eps}(\tau)}{h}\\ \notag
=& \ \Big(\frac{p_{t_0,\tau+h}^w(x,D_x)-p_{t_0,\tau}^w(x,D_x)}{h}\chi_{\eps}^w(x,D_x)\mathscr{U}(\tau+h,\tau_0)u,v\Big)_{L^2(\rr^n)}\\ \notag
& \ + \Big(p_{t_0,\tau}^w(x,D_x)\chi_{\eps}^w(x,D_x)\frac{\mathscr{U}(\tau+h,\tau_0)-\mathscr{U}(\tau,\tau_0)}{h}u,v\Big)_{L^2(\rr^n)}.
\end{align}
By using anew that the mappings $\chi_{\eps}^w(x,D_x) : L^2(\rr^n) \rightarrow \mathscr{S}(\rr^n)$ and $p_{t_0,\tau}^w(x,D_x) : \mathscr{S}(\rr^n) \rightarrow \mathscr{S}(\rr^n)$ are continuous thanks to Proposition~\ref{FIOdef} and (\ref{jk31}), we deduce from Definition~\ref{def1} and Theorem~\ref{th1} that 
\begin{multline}\label{as7}
\lim_{h \to 0} \Big(p_{t_0,\tau}^w(x,D_x)\chi_{\eps}^w(x,D_x)\frac{\mathscr{U}(\tau+h,\tau_0)-\mathscr{U}(\tau,\tau_0)}{h}u,v\Big)_{L^2(\rr^n)}\\
=\big(p_{t_0,\tau}^w(x,D_x)\chi_{\eps}^w(x,D_x)q_\tau^w(x,D_x)\mathscr{U}(\tau,\tau_0)u,v\big)_{L^2(\rr^n)},
\end{multline}
since $\tau \mapsto \mathscr{U}(\tau,\tau_0)u \in C^1(]\tau_0,t_0],L^2(\rr^n))$.
On the other hand, it follows from (\ref{wigner}) and (\ref{hj2}) that 
\begin{multline*}
\Big(\frac{p_{t_0,\tau+h}^w(x,D_x)-p_{t_0,\tau}^w(x,D_x)}{h}\chi_{\eps}^w(x,D_x)\mathscr{U}(\tau+h,\tau_0)u,v\Big)_{L^2(\rr^n)}\\
=\frac{1}{h}\int_{\rr^{2n}}\big(p_{t_0,\tau+h}(x,\xi)-p_{t_0,\tau}(x,\xi)\big)
\mathcal{H}\big(\chi_{\eps}^w(x,D_x)\mathscr{U}(\tau+h,\tau_0)u,v\big)(x,\xi)dxd\xi,
\end{multline*}
since $p_{t,\tau}$ is a $L^{\infty}(\rr^{2n})$-function when $0 \leq \tau \leq t \leq T$, $0 \leq t-\tau \leq \delta$.
The above integral is well-defined as the Wigner function $\mathcal{H}\big(\chi_{\eps}^w(x,D_x)\mathscr{U}(\tau+h,\tau_0)u,v\big)$ belongs to the Schwartz space $\mathscr{S}(\rr^{2n})$ since $\chi_{\eps}^w(x,D_x)\mathscr{U}(\tau+h,\tau_0)u \in \mathscr{S}(\rr^{n})$ and $v \in \mathscr{S}(\rr^{n})$. The continuity of the mapping $h \mapsto \mathscr{U}(\tau+h,\tau_0)u \in L^2(\rr^n)$ successively implies the continuity of the mappings $h \mapsto \chi_{\eps}^w(x,D_x)\mathscr{U}(\tau+h,\tau_0)u \in \mathscr{S}(\rr^n)$ and $h \mapsto \mathcal{H}\big(\chi_{\eps}^w(x,D_x)\mathscr{U}(\tau+h,\tau_0)u,v\big) \in \mathscr{S}(\rr^{2n})$. We therefore deduce from (\ref{carre}) and (\ref{ku1}) that 
\begin{align}\label{as8}
& \ \lim_{h \to 0}\Big(\frac{p_{t_0,\tau+h}^w(x,D_x)-p_{t_0,\tau}^w(x,D_x)}{h}\chi_{\eps}^w(x,D_x)\mathscr{U}(\tau+h,\tau_0)u,v\Big)_{L^2(\rr^n)}\\ \notag
=& \ -\int_{\rr^{2n}}\big(p_{t_0,\tau} \#^w q_\tau)(x,\xi)\mathcal{H}\big(\chi_{\eps}^w(x,D_x)\mathscr{U}(\tau,\tau_0)u,v\big)(x,\xi)dxd\xi\\ \notag
=& \ -\big(p_{t_0,\tau}^w(x,D_x)q_\tau^w(x,D_x)\chi_{\eps}^w(x,D_x)\mathscr{U}(\tau,\tau_0)u,v\big)_{L^2(\rr^n)},
\end{align}
since the domination condition follows as above from the fact that any Schwartz seminorm 
of the Wigner function can be bounded as 
\begin{multline*}
\sup_{\substack{x,\xi \in \rr^n, \\ |\alpha_1|+|\alpha_2|+|\beta_1|+|\beta_2| \leq N_1}}|x^{\alpha_1}\xi^{\alpha_2}\partial_{x}^{\beta_1}\partial_{\xi}^{\beta_2}\mathcal{H}(\chi_{\eps}^w(x,D_x)\mathscr{U}(\tau+h,\tau_0)u,v)(x,\xi)| \\
\leq c \|\mathscr{U}(\tau+h,\tau_0)u\|_{L^2(\rr^n)}\Big(\sup_{\substack{x \in \rr^n, \\ |\alpha|+|\beta| \leq N_2}}|x^{\alpha}\partial_x^{\beta}v(x)|\Big)
\leq c \|u\|_{L^2(\rr^n)}\Big(\sup_{\substack{x \in \rr^n, \\ |\alpha|+|\beta| \leq N_2}}|x^{\alpha}\partial_x^{\beta}v(x)|\Big),
\end{multline*}
since $\|\mathscr{U}(\tau+h,\tau_0)\|_{\mathcal{L}(L^2)} \leq 1$.
It follows from (\ref{as9}), (\ref{as7}) and (\ref{as8}) that for all $\tau_0<\tau<t_0$,
\begin{equation}\label{as10}
f_{\eps}'(\tau)=\big(p_{t_0,\tau}^w(x,D_x)[\chi_{\eps}^w(x,D_x),q_\tau^w(x,D_x)]\mathscr{U}(\tau,\tau_0)u,v\big)_{L^2(\rr^n)}.
\end{equation}
We deduce from (\ref{carre}), (\ref{as1}) and (\ref{as10}) that 
\begin{multline}\label{as11}
\big(\chi_{\eps}^w(x,D_x)\mathscr{U}(t_0,\tau_0)u,v\big)_{L^2(\rr^n)}-\big(p_{t_0,\tau_0}^w(x,D_x)\chi_{\eps}^w(x,D_x)u,v\big)_{L^2(\rr^n)}\\
=\int_{\tau_0}^{t_0}\big(p_{t_0,\tau}^w(x,D_x)[\chi_{\eps}^w(x,D_x),q_\tau^w(x,D_x)]\mathscr{U}(\tau,\tau_0)u,v\big)_{L^2(\rr^n)}d\tau,
\end{multline}
since $\mathscr{U}(\tau_0,\tau_0)=I_{L^2(\rr^n)}$.
By passing to the limit when $\eps$ tends to $0$, it follows from Proposition~\ref{FIOdef}, (\ref{jk31}) and (\ref{tria}) that 
\begin{multline}\label{as12}
\big(\big(\mathscr{U}(t_0,\tau_0)-p_{t_0,\tau_0}^w(x,D_x)\big)u,v\big)_{L^2(\rr^n)}\\
=\lim_{\eps \to 0}\int_{\tau_0}^{t_0}\big(p_{t_0,\tau}^w(x,D_x)[\chi_{\eps}^w(x,D_x),q_\tau^w(x,D_x)]\mathscr{U}(\tau,\tau_0)u,v\big)_{L^2(\rr^n)}d\tau,
\end{multline}
since $\mathscr{U}(t_0,\tau_0)$ and $p_{t_0,\tau_0}^w(x,D_x)$ are bounded operators on $L^2(\rr^n)$. By using that the Weyl symbol of the operator $q_{\tau}^w(x,D_x)$ is quadratic and (\ref{carre2}), standard results of symbolic calculus show that the commutator $[\chi_{\eps}^w(x,D_x),q_{\tau}^w(x,D_x)]$ is equal to 
\begin{multline}\label{as13}
[\chi_{\eps}^w(x,D_x),q_\tau^w(x,D_x)]=\sum_{\substack{\alpha, \beta \in \mathbb{N}^n \\ |\alpha+\beta| = 2}}(q_\tau)_{\alpha,\beta}[\chi_{\eps}^w(x,D_x),(x^{\alpha}\xi^{\beta})^w]\\
=\frac{1}{i}\sum_{\substack{\alpha, \beta \in \mathbb{N}^n \\ |\alpha+\beta| = 2}}(q_\tau)_{\alpha,\beta}\textrm{Op}^w\big(\{\chi_{\eps},x^{\alpha}\xi^{\beta}\}\big).
\end{multline}
We notice that the symbol
$$\{\chi_{\eps},x^{\alpha}\xi^{\beta}\}(x,\xi)=\eps \sum_{j=1}^n\Big(\frac{\partial \chi}{\partial \xi_j}(\eps x,\eps \xi) \cdot \frac{\partial (x^{\alpha}\xi^{\beta})}{\partial x_j}-\frac{\partial \chi}{\partial x_j}(\eps x,\eps \xi) \cdot \frac{\partial (x^{\alpha}\xi^{\beta})}{\partial \xi_j}\Big),$$
writes as $\Psi_{\eps}(x,\xi)=\Psi(\eps x,\eps \xi)$, with $\Psi \in C_0^{\infty}(\rr^{2n},\cc)$.
It is therefore uniformly bounded in the Fr\'echet space $C_b^{\infty}(\rr^{2n})$ with respect to $0<\eps \leq 1$. On the other hand, this symbol vanishes on any compact set when $0<\eps \ll 1$. It therefore converges in the Fr\'echet space $C^{\infty}(\rr^{2n})$ to zero when $\eps$ tends to $0$. By using the very same arguments as in (\ref{conv}), we obtain that 
\begin{equation}\label{as20}
\forall w \in L^2(\rr^n), \quad \lim_{\eps \to 0}\|\textrm{Op}^w\big(\{\chi_{\eps},x^{\alpha}\xi^{\beta}\}\big)w\|_{L^2(\rr^n)}=0.
\end{equation}
Furthermore, the Calder\'on-Vaillancourt Theorem together with the continuity of the coefficients $\tau \in [0,T] \mapsto (q_{\tau})_{\alpha,\beta} \in \cc$ imply that there exists a positive constant $C_1>0$ such that 
\begin{equation}\label{as30}
\forall \tau_0 \leq \tau \leq t_0, \forall 0<\eps \leq 1, \quad \|[\chi_{\eps}^w(x,D_x),q_\tau^w(x,D_x)]\|_{\mathcal{L}(L^2(\rr^n))} \leq C_1.
\end{equation}
Recalling from Proposition~\ref{FIOdef} and (\ref{jk31}) that $p_{t_0,\tau}^w(x,D_x)$ defines a bounded operator on $L^2(\rr^n)$, 
we deduce from (\ref{as13}) and (\ref{as20}) that for all $\tau_0 \leq \tau \leq t_0$,   
\begin{equation}\label{as40}
\lim_{\eps \to 0}\big(p_{t_0,\tau}^w(x,D_x)[\chi_{\eps}^w(x,D_x),q_\tau^w(x,D_x)]\mathscr{U}(\tau,\tau_0)u,v\big)_{L^2(\rr^n)}=0.
\end{equation}
On the other hand, it follows from (\ref{as30}) that for all $\tau_0 \leq \tau \leq t_0$,
\begin{equation}\label{as41}
\big|\big(p_{t_0,\tau}^w(x,D_x)[\chi_{\eps}^w(x,D_x),q_\tau^w(x,D_x)]\mathscr{U}(\tau,\tau_0)u,v\big)_{L^2(\rr^n)}\big| \leq C_1\|u\|_{L^2(\rr^n)}\|v\|_{L^2(\rr^n)},
\end{equation}
since from Theorem~\ref{th1}, Proposition~\ref{FIOdef} and (\ref{jk31}), we have $\|p_{t_0,\tau}^w(x,D_x)\|_{\mathcal{L}(L^2)} \leq 1$ and $\|\mathscr{U}(\tau,\tau_0)\|_{\mathcal{L}(L^2)} \leq 1$.
By Lebesgue's theorem, we deduce from (\ref{as12}), (\ref{as40}) and (\ref{as41}) that 
\begin{equation}\label{as42}
\forall u, v \in \mathscr{S}(\rr^n), \quad \big(\big(\mathscr{U}(t_0,\tau_0)-p_{t_0,\tau_0}^w(x,D_x)\big)u,v\big)_{L^2(\rr^n)}=0.
\end{equation}
By density of the Schwartz space in $L^2(\rr^n)$ and the continuity of the operators on $L^2(\rr^n)$, we finally conclude that 
\begin{equation}\label{as43}
\forall u\in L^2(\rr^n), \quad \mathscr{U}(t_0,\tau_0)u=p_{t_0,\tau_0}^w(x,D_x)u,
\end{equation}
that is, $\mathscr{U}(t_0,\tau_0)=p_{t_0,\tau_0}^w(x,D_x)$. This ends the proof of Theorem~\ref{th3}.

On the other hand, let $0 \leq \tau \leq t \leq T$. We choose a finite sequence $(s_j)_{1 \leq j \leq N}$, with $N \geq 2$ satisfying 
$$s_1=\tau < s_2 <....< s_{N-1} < s_N=t, \qquad  0 < s_{j+1}-s_j < \delta, \quad 1 \leq j \leq N-1,$$
where $\delta>0$ is the positive constant given by Theorem~\ref{th3}.
We deduce from (\ref{jk31}), Theorems~\ref{th1} and~\ref{th3} that 
\begin{multline}\label{tri10}
\mathscr{U}(t,\tau)=\mathscr{U}(s_N,s_1)=\mathscr{U}(s_N,s_{N-1})...\mathscr{U}(s_2,s_1)\\
=p_{s_N,s_{N-1}}^w(x,D_x)...p_{s_2,s_1}^w(x,D_x)=\mathscr{K}_{R(s_N,s_{N-1})}...\mathscr{K}_{R(s_2,s_1)}.
\end{multline}
It is shown in~\cite{mehler} (Proposition~5.9) that if $\mathcal{T}_1$ and $\mathcal{T}_2$ are non-negative complex symplectic linear transformations then $\mathcal{T}_1\mathcal{T}_2$ is also a non-negative complex symplectic linear transformation and the associated Fourier integral operators satisfy either
$$\mathscr{K}_{\mathcal{T}_1\mathcal{T}_2}=\mathscr{K}_{\mathcal{T}_1}\mathscr{K}_{\mathcal{T}_2}$$
or 
$$\mathscr{K}_{\mathcal{T}_1\mathcal{T}_2}=-\mathscr{K}_{\mathcal{T}_1}\mathscr{K}_{\mathcal{T}_2}.$$
Recalling from Proposition~\ref{FIOdef} that the kernels of the Fourier integral operators are only determined up to their signs, we may therefore consider that the following formula holds true
\begin{equation}\label{tri11}
\mathscr{K}_{\mathcal{T}_1\mathcal{T}_2}=\mathscr{K}_{\mathcal{T}_1}\mathscr{K}_{\mathcal{T}_2},
\end{equation}
whenever $\mathcal{T}_1$ and $\mathcal{T}_2$ are non-negative complex symplectic linear transformations. We therefore deduce from (\ref{tri10}) and (\ref{tri11}) that 
$$\mathscr{U}(t,\tau)=\mathscr{K}_{R(s_N,s_{N-1})}...\mathscr{K}_{R(s_2,s_1)}=\mathscr{K}_{R(s_N,s_{1})}=\mathscr{K}_{R(t,\tau)}.$$
Theorem~\ref{th2} then directly follows from Proposition~\ref{FIOdef}.

\section{Propagation of Gabor singularities}\label{gabor}

This section is devoted to give the proof of Theorem~\ref{th4}. Let $T>0$ and $q_t : \rr^{2n} \rightarrow \cc$ be a time-dependent complex-valued quadratic form with a non-positive real part $\textrm{Re }q_t \leq 0$ for all 
$0 \leq t \leq T$, and
whose coefficients depend continuously on the time variable $0 \leq t \leq T$.

We aim at studying the possible (or lack of) Schwartz regularity for the $B$-valued solutions $u(t)=\mathscr{U}(t,0)u_0$ at time $0 \leq t \leq T$ to the non-autonomous Cauchy problem 
$$\left\lbrace\begin{array}{l}
\frac{du(t)}{dt}=q_t^w(x,D_x)u(t), \qquad 0 < t \leq T,\\
u(0)=u_0,
\end{array}\right.$$
given by Theorem~\ref{th1}, where $u_0 \in B$ is an arbitrary initial datum. To that end, we derive a microlocal inclusion for the Gabor wave front set of the solution $u(t)=\mathscr{U}(t,0)u_0$ in terms of the Hamilton maps $(F_{\tau})_{0 \leq \tau \leq t}$ of the quadratic symbols $(q_{\tau})_{0 \leq \tau \leq t}$ and the Gabor wave front set of the initial datum $WF(u_0)$.  Thanks to Theorem~\ref{th2}, the proof of Theorem~\ref{th4} is an adaptation of the analysis led in~\cite{wahlberg} in the autonomous case. The keystone in~\cite{wahlberg} (Theorem~4.6) is the proof of the microlocal inclusion 
\begin{equation}\label{incl1}
WF(K_{\mathcal{T}}) \subset (\widetilde{\lambda_{\mathcal{T}}} \cap \rr^{4n}) \setminus \{0\},
\end{equation}
for the Gabor wave front set of $K_{\mathcal{T}} \in \mathscr{S}'(\rr^{2n})$ the kernel of the Fourier integral operator $\mathscr{K}_{\mathcal{T}}$ defined in Proposition~\ref{FIOdef} and associated to a non-negative complex symplectic linear transformation $\mathcal{T}$, where $\widetilde{\lambda_{\mathcal{T}}}$ denotes the non-negative Lagrangian plane (\ref{graph1}). It follows from (\ref{graph1}) and (\ref{incl1}) that  
\begin{multline}\label{incl2}
WF(K_{\mathcal{T}}) \\
\subset \big\{(x,y,\xi,-\eta) \in \rr^{4n} \setminus \{0\} : (x,\xi)=\mathcal{T}(y,\eta), \ (y,\eta) \in \textrm{Ker}(\textrm{Im }\mathcal{T}) \cap \rr^{2n} \big\},
\end{multline}
with $\textrm{Im }\mathcal{T}=\frac{1}{2i}(\mathcal{T}-\overline{\mathcal{T}})$. We notice from (\ref{incl2}) that the Gabor wave front set of the kernel $K_{\mathcal{T}} \in \mathscr{S}'(\rr^{2n})$ does not contain any point of the form $(0,y,0,-\eta)$ for $(y,\eta) \in \rr^{2n} \setminus \{0\}$, nor points of the form $(x,0,\xi,0)$ for $(x,\xi) \in \rr^{2n} \setminus \{0\}$, since $\mathcal{T} : \cc^{2n} \rightarrow \cc^{2n}$ is invertible. We can therefore deduce from~\cite{Hormander1} (Proposition~2.11) the microlocal inclusion
\begin{equation}\label{composition1}
\forall u \in \mathscr{S}'(\rr^n), \quad WF(\mathscr{K}_{\mathcal{T}}u) 
\subset  WF'(K_{\mathcal{T}}) \circ WF(u), 
\end{equation}
that is,
\begin{multline}\label{composition1bis}
\forall u \in \mathscr{S}'(\rr^n), \quad WF(\mathscr{K}_{\mathcal{T}}u) \\
\subset  \big\{(x,\xi)\in \rr^{2n} \setminus \{0\} : \exists (y,\eta) \in WF(u), \ (x,y,\xi,-\eta) \in WF(K_{\mathcal{T}})\big\}.
\end{multline}
It follows from (\ref{incl2}) and (\ref{composition1bis}) that 
\begin{multline*}
\forall u \in \mathscr{S}'(\rr^n), \quad WF(\mathscr{K}_{\mathcal{T}}u) 
\subset \big\{(x,\xi) \in \rr^{2n} \setminus \{0\} : \\
\exists (y,\eta) \in WF(u) \cap \textrm{Ker}(\textrm{Im }\mathcal{T}) \cap \rr^{2n}, \  (x,\xi)=\mathcal{T}(y,\eta)\big\}, 
\end{multline*}
that is
\begin{equation}\label{composition3}
\forall u \in \mathscr{S}'(\rr^n), \quad WF(\mathscr{K}_{\mathcal{T}}u) 
\subset \mathcal{T}\big( WF(u) \cap \textrm{Ker}(\textrm{Im }\mathcal{T}) \cap \rr^{2n}\big).
\end{equation}
By noticing that 
\begin{equation}\label{ku6}
\mathcal{T}\big( WF(u) \cap \textrm{Ker}(\textrm{Im }\mathcal{T}) \cap \rr^{2n}\big)=\mathcal{T}\big(WF(u)\big) \cap \textrm{Ker}(\textrm{Im }\mathcal{T}^{-1}) \cap \rr^{2n},
\end{equation}
it follows from (\ref{composition3}) and (\ref{ku6}) that 
\begin{equation}\label{ku6.5}
\forall u \in \mathscr{S}'(\rr^n), \quad WF(\mathscr{K}_{\mathcal{T}}u) 
\subset \mathcal{T}\big(WF(u)\big) \cap \textrm{Ker}(\textrm{Im }\mathcal{T}^{-1}) \cap \rr^{2n}.
\end{equation}
On the other hand, we deduce from (\ref{ku6.5}), Theorems~\ref{th1} and~\ref{th2} that for all $u_0 \in B$ and $0 \leq \tau \leq  t \leq T$, 
\begin{align}\label{ku5}
& \ WF(\mathscr{U}(t,0)u_0)=WF(\mathscr{U}(t,\tau)\mathscr{U}(\tau,0)u_0)\\ \notag
\subset & \ R(t,\tau)\big[WF(\mathscr{U}(\tau,0)u_0)\big]  \cap \textrm{Ker}\big(\textrm{Im }R(\tau,t)\big) \cap \rr^{2n}\\ \notag
\subset & \ R(t,\tau)\Big[R(\tau,0)\big( WF(u_0)\big)  \cap \textrm{Ker}\big(\textrm{Im }R(0,\tau)\big) \cap \rr^{2n}\Big]  \cap \textrm{Ker}(\textrm{Im }R(\tau,t)) \cap \rr^{2n}\\ \notag
\subset & \ R(t,0)\big( WF(u_0)\big)  \cap \textrm{Ker}(\textrm{Im }R(\tau,t)) \cap \rr^{2n}.
\end{align}
Then, it follows from (\ref{ku5}) that for all $0 \leq t \leq T$,
\begin{equation}\label{ku7}
WF(\mathscr{U}(t,0)u_0) \subset  R(t,0)\big(WF(u_0)\big)  \cap S_{0,t}.
\end{equation}
where $S_{0,t}$ is the time-dependent singular space
\begin{equation}\label{ku8}
S_{0, t}=\Big( \bigcap_{0 \leq \tau \leq t}\textrm{Ker}(\textrm{Im }R(\tau,t))\Big) \cap \rr^{2n},
\end{equation}
defined in Definition~\ref{def2}. With 
$$\textrm{Re }R(t,0)=\frac{1}{2}(R(t,0)+\overline{R(t,0)}),$$ 
we finally obtain that for all $0 \leq t \leq T$ and $u_0 \in B$, 
\begin{equation}\label{ku9}
WF(u(t))=WF(\mathscr{U}(t,0)u_0) \subset  \big(\textrm{Re }R(t,0)\big)\big(WF(u_0)\big)  \cap S_{0,t},
\end{equation}
since $WF(u_0) \subset \rr^{2n} \setminus \{0\}$ and $S_{0,t} \subset \rr^{2n}$. This ends the proof of Theorem~\ref{th4}.

\section{Appendix: Gabor wave front set}\label{appendix}

This appendix is devoted to recall the definition and basic properties of the Gabor wave front set of a tempered distribution. This wave front set is defined as a subset of the phase space characterizing the lack of Schwartz regularity of the tempered distribution. 

For all $x,y,\xi \in \rr^n$, we denote
$$T_x f(y)=f(y-x), \quad M_\xi f(y)=e^{i  y \cdot \xi } f(y), \quad \Pi(x,\xi) = M_\xi T_x,$$
the translation, modulation and phase space translation operators. Given a window function $\varphi \in \mathscr{S}(\rr^n) \setminus \{0\}$, the short-time Fourier transform of the tempered distribution $f \in \mathscr{S}'(\rr^n)$ is defined in~\cite{Grochenig1} as 
$$(V_\varphi f)(x,\xi) = \langle f, \overline{\Pi(x,\xi) \varphi} \rangle_{\mathscr{S}'(\rr^n),\mathscr{S}(\rr^n)}, \qquad (x,\xi) \in \rr^{2n}.$$
The function $(x,\xi) \in \rr^{2n} \mapsto (V_\varphi f)(x,\xi) \in \cc$ is smooth and its modulus is bounded by $C \langle (x,\xi) \rangle^k$ for all $(x,\xi) \in \rr^{2n}$ for some constants $C,k \geq 0$.
If $\varphi \in \mathscr{S}(\rr^n)$, $\|\varphi\|_{L^2(\rr^n)}=1$ and $f \in \mathscr{S}'(\rr^n)$,
the short-time Fourier transform inversion formula~\cite[Corollary 11.2.7]{Grochenig1} reads as
$$\forall g \in \mathscr{S}(\rr^n), \quad \langle f,g\rangle_{\mathscr{S}'(\rr^n),\mathscr{S}(\rr^n)}  = \frac{1}{(2 \pi)^{n}} \int_{\rr^{2n}} (V_\varphi f)(x,\xi) \langle \Pi(x,\xi) \varphi,g\rangle_{\mathscr{S}'(\rr^n),\mathscr{S}(\rr^n)} dxd\xi.$$
On the other hand, we recall that the Shubin symbol class $G^m$, with $m \in \rr$, is defined as the space of all $a \in C^\infty(\rr^{2n},\cc)$ satisfying 
\begin{equation}\label{symbolestimate1}
\forall \alpha,\beta \in \mathbb{N}^n, \exists C_{\alpha,\beta}>0, \forall (x,\xi)\in \rr^{2n}, \quad  |\partial_x^\alpha \partial_\xi^\beta a(x,\xi)| 
\leq C_{\alpha,\beta} \langle (x,\xi) \rangle^{m-|\alpha|-|\beta|}. 
\end{equation}
The space $G^m$ equipped with the semi-norms 
$$\sup_{(x,\xi) \in \rr^{2n}} \langle (x,\xi) \rangle^{-m+|\alpha|+|\beta|}| \partial_x^\alpha \partial_\xi^\beta a(x,\xi)|, \quad \alpha,\beta \in \mathbb{N}^{n},$$
is a Fr\'echet space. Given a Shubin symbol $a \in G^m$, a non-zero point in the phase space $(x_0,\xi_0) \in \rr^{2n} \setminus \{ (0,0) \}$ is said to be non-characteristic for the symbol $a$ with respect to the class $G^m$ provided there exist some positive constants $A,\eps>0$ and an open conic\footnote{A set invariant under multiplication with positive reals.} set $\Gamma \subseteq \rr^{2n} \setminus \{ (0,0) \}$ such that 
$$(x_0,\xi_0) \in \Gamma, \qquad  \forall (x,\xi) \in \Gamma, \forall |(x,\xi)| \geq A, \quad  |a(x,\xi)| \geq \eps \langle (x,\xi) \rangle^m.$$
Otherwise, the non-zero point $(x_0,\xi_0) \in \rr^{2n} \setminus \{ (0,0) \}$ is said to be characteristic. We denote by $\textrm{Char}(a) \subset \rr^{2n} \setminus \{ (0,0) \}$ the set of all characteristic points.

The notion of Gabor wave front set is defined as follows by H\"ormander~\cite{Hormander1} to measure the directions in the phase space in which a tempered distribution does not behave like a Schwartz function: 

\bigskip

\begin{definition}\label{wavefront1}
Let $u \in \mathscr{S}'(\rr^n)$ be a tempered distribution. Its Gabor wave front set $WF(u)$ is defined as the set of all non-zero points in the phase space $(x,\xi) \in \rr^{2n} \setminus \{ (0,0) \}$ such that for all $a \in G^m$, with $m \in \rr$,
$$a^w(x,D_x) u \in \mathscr{S}(\rr^n)   \Longrightarrow (x,\xi) \in \emph{\textrm{Char}}(a).$$
\end{definition}

\bigskip

According to \cite[Proposition 6.8]{Hormander1} and \cite[Corollary 4.3]{Rodino1}, the Gabor wave front set can be microlocally characterized by the short-time Fourier transform. Indeed, if $u \in \mathscr{S}'(\rr^n)$ and $\varphi \in \mathscr{S}(\rr^n) \setminus \{0\}$, then $(x_0,\xi_0) \in \rr^{2n} \setminus \{(0,0)\}$ satisfies $(x_0,\xi_0) \notin WF(u)$ if and only if there exists an open conic set $\Gamma_{x_0,\xi_0} \subseteq \rr^{2n} \setminus \{ (0,0) \}$ containing $(x_0,\xi_0)$ such that
$$\forall N \geq 0, \quad \sup_{(x,\xi) \in \Gamma_{x_0,\xi_0}} \langle (x,\xi) \rangle^N |(V_\varphi u)(x,\xi)| < +\infty.$$
The Gabor wave front set satisfies the following basic properties: 
\begin{enumerate}
\item[$(i)$] If $u \in \mathscr{S}'(\rr^n)$, then \cite[Proposition 2.4]{Hormander1}
\begin{equation}\label{nm1}
WF(u) = \emptyset \Longleftrightarrow  u \in \mathscr{S}(\rr^n)
\end{equation} 
\item[$(ii)$] If $u \in \mathscr{S}'(\rr^n)$ and $a \in G^m$, then
$$WF(a^w(x,D_x) u) 
 \subset WF(u) \cap \textrm{conesupp}(a) 
 \subset WF( a^w(x,D_x) u) \ \cup \ \textrm{Char}(a),$$
where the conic support $\textrm{conesupp}(a)$ of $a \in G^m$  is the set of all
$(x,\xi) \in \rr^{2n} \setminus \{0\}$ such that any conic open set $\Gamma_{x,\xi} \subseteq \rr^{2n} \setminus \{0\}$ containing $(x,\xi)$ verifies
$$\overline{\textrm{supp}(a) \cap \Gamma_{x,\xi}} \quad \mbox{is not compact in} \ \rr^{2n}$$
\end{enumerate}
The Gabor wave front set also enjoys some symplectic invariant features thanks to the symplectic invariance of the Weyl quantization.
We recall that the real symplectic group $\textrm{Sp}(n,\rr)$ consists of all matrices $\chi \in \textrm{GL}(2n,\rr)$ preserving the symplectic form
\begin{equation}\label{symplecticgroup}
\sigma\big(\chi(X),\chi(X')\big) = \sigma(X,X'),
\end{equation}
for all $X,X' \in \rr^{2n}$, whereas the complex symplectic group $\textrm{Sp}(n,\cc)$ consists of all matrices $\chi \in \textrm{GL}(2n,\cc)$ satisfying \eqref{symplecticgroup} for all $X,X' \in \cc^{2n}$. 
To each real symplectic matrix $\chi \in \textrm{Sp}(n,\rr)$ is associated~\cite{Folland1,Hormander0} a unitary operator $\mu(\chi)$ on $L^2(\rr^n)$, determined up to a complex factor of modulus one, satisfying
\begin{equation}\label{symplecticoperator}
\forall a \in \mathscr{S}'(\rr^{2n}), \quad \mu(\chi)^{-1} a^w(x,D_x)  \mu(\chi) = (a \circ \chi)^w(x,D_x).
\end{equation}
The operator $\mu(\chi)$ is an homeomorphism on $\mathscr{S}(\rr^n)$ and on $\mathscr{S}'(\rr^n)$.
The mapping $\textrm{Sp}(n,\rr) \ni \chi \mapsto \mu(\chi)$ is called the metaplectic representation~\cite{Folland1}.
It is in fact a representation of the so called $2$-fold covering group of $\textrm{Sp}(n,\rr)$, which is called the metaplectic group and denoted $\textrm{Mp}(n,\rr)$.
The metaplectic representation satisfies the homomorphism relation only modulo a change of sign
$$\mu( \chi \chi') = \pm \mu(\chi ) \mu(\chi' ), \quad \chi, \chi' \in \textrm{Sp}(n,\rr).$$
According to \cite[Proposition~2.2]{Hormander1}, the Gabor wave front set is symplectically invariant, that is, for all $u \in \mathscr{S}'(\rr^n)$, $\chi \in \textrm{Sp}(n,\rr)$,
$$(x,\xi) \in WF(u) \quad \Longleftrightarrow \quad \chi(x,\xi) \in WF(\mu(\chi)u),$$
that is,
\begin{equation}\label{sympinv}
WF( \mu(\chi) u) = \chi WF(u), \quad \chi \in \textrm{Sp}(n,\rr), \ u \in \mathscr{S}'(\rr^n).
\end{equation}

\end{document}